\documentclass[final]{siamltex}
\usepackage{animate}
\usepackage{amsmath}
\usepackage{amssymb}
\usepackage{graphics}
\usepackage{graphicx}
\usepackage{footnote}
\usepackage{textcomp}
\usepackage{mathrsfs}
\usepackage{epstopdf}
\usepackage{array}
\usepackage{showkeys}
\usepackage[maxfloats=99]{morefloats}
\usepackage{url}
\usepackage{cases}
\usepackage{mathscinet}
\usepackage[normalem]{ulem}
\usepackage{algorithmic, algorithm}
\usepackage{color}
\usepackage{cite}

\usepackage{subcaption}

\newcommand{\R}{\mathbb R}
\newcommand{\N}{\mathbb N}

\newcommand{\f}[1]{\mathbf{#1}}
\newcommand{\norm}[1]{\left\| #1 \right\|}
\newcommand{\abs}[1]{\lvert#1\rvert}
\newcommand{\set}[1]{\left\{ #1\right\}}
\newcommand{\lapla}{\Delta}
\newcommand{\manifold}{\mathcal{M}}
\newcommand{\I}{\text{Id}}
\numberwithin{equation}{section}

\newtheorem{assumption}[theorem]{Assumption}

\newtheorem{remark}{Remark}[section]

\usepackage[T1]{fontenc}
\usepackage{textcomp}
\usepackage[scaled]{helvet}

\title{Parametric Polynomial Preserving Recovery on Manifolds }
\author{Guozhi Dong\thanks{Computational Science Center, Faculty of Mathematics, University of Vienna, 
Oskar-Morgenstern-Platz 1, 1090 Vienna, Austria,
and Institute for Mathematics, Humboldt University of Berlin, Unter den Linden 6, 10099 Berlin, Germany (guozhi.dong@hu-berlin.de).  This author was partially supported by  the Austrian Science Fund (FWF): Geometry and Simulation, project S11704} 
\and %
Hailong Guo\thanks{School of Mathematics and Statistics, The University of Melbourne, Parkville, VIC, 3010, Australia (hailong.guo@unimelb.edu.au). This author  was partially supported by Andrew Sisson Fund of the University of Melbourne.
} 
}

\begin{document}
\maketitle
\medskip

\begin{abstract}
This paper investigates gradient recovery schemes for data defined on discretized manifolds.
The proposed method, parametric polynomial preserving recovery (PPPR), does not require the tangent spaces of the exact manifolds, and they have been assumed for some significant gradient recovery methods in the literature.
Another advantage of PPPR is that superconvergence is guaranteed without the symmetric condition which has been asked in the existing techniques.
There is also numerical evidence that the superconvergence by PPPR is high curvature stable, which distinguishes itself from the others.
As an application, we show its capability of constructing an asymptotically exact \textit{a posteriori} error estimator.
Several numerical examples on two-dimensional surfaces are presented to support the theoretical results and comparisons with  existing methods are documented.
\vskip .3cm
{\bf AMS subject classifications.} \ {Primary 65N50, 65N30; Secondary 65N15, 53C99}
\vskip .3cm

{\bf Key words.} \ {Gradient recovery, manifolds, superconvergence, parametric polynomial preserving, function value preserving, curvature stable.}
\end{abstract}


\section{Introduction}
Numerical methods for approximating variational problems or partial differential equations (PDEs) with solutions defined on surfaces or manifolds are of growing interests over the last decades.
Finite element methods, as one of the most important methods for numerically solving PDEs, are well established for those problems.
A starting point can be traced back to \cite{Dziuk1988}, which is the first to investigate a finite element method for solving elliptic PDEs on surfaces.
Since then, there have been a lot of extensions in both numerical analysis and practical algorithms,  see \cite{DemlowDziuk2007,Demlow2009,OlshanskiiReuskenGrande2009,OlshankiiReusken2010,DednerMadhavanStinner2013,DziukElliott2013,OlshanskiiSafin2016} and the references therein.
In the literature, most  works focus on the {\it  a priori} error analysis of various surface finite element methods. Only a few works, up to our best knowledge, take into account the {\it a posteriori}  error analysis and superconvergence of finite element methods in a surface setting, see \cite{DuJu2005,DemlowDziuk2007,WeiChenHuang2010,Demlow2012,CamachoDemlow2015,Chernyshenk02015,DednerMadhavanStinner2016}.
Recently, there is an approach proposed in \cite{GrandeReusken2016} which merges the two types of analysis to develop a higher order finite element method on an approximated surface, where a gradient recovery scheme plays a vital role.
Gradient recovery techniques, which are important in \emph{post-processing} solutions or data to improve the accuracy of gradient approximations, have been widely studied and found many applications in numerical analysis.
For planar problems, the study of gradient recovery methods has reached a certain maturity stage,
and there is a massive of works in the literature, to name  a few \cite{BankXu2003,XuZhang2004,AinsworthOden2000,GuoZhang2015,Lakhany2000,ZhangNaga2005,ZZ1992,ZZ1992b}.
We point out some significant methods among them, like the classical Zienkiewicz--Zhu ($ZZ$) superconvergent patch recovery  \cite{ZZ1992},
and a later method called polynomial preserving recovery (PPR) \cite{ZhangNaga2005}.
Those two approaches work under different philosophies methodologically.
The former method first locates positions of superconvergent  points in the neighborhood of each nodal point, and then recovers the gradients themselves at the nodal point by fitting the selected neighbored superconvergent points to achieve a higher order approximation accuracy;
while the latter first fits a polynomial in the least-squares sense at each nodal point and then takes the gradient of the fitted polynomial to have the recovered gradient. 
Both methods can produce comparable superconvergence results, but $ZZ$-Scheme requires stronger conditions on the discretized meshes than the PPR method.

Gradient recovery methods for data defined on curved spaces have only recently been investigated.
In \cite{WeiChenHuang2010}, several gradient recovery methods have been extended to a general surface setting for linear finite element solutions which are defined on polyhedrons by triangulation.
The surface in \cite{WeiChenHuang2010} is considered to be a zero level set of a smooth function embedded in an ambient space, and the gradient of functions on the surface is represented using the ambient gradient with tangential projection.
It has been shown that most of the properties of the gradient recovery schemes for planar problems are maintained in their counterparts for surface problems.
In particular, in their implementation and analysis, the methods require knowledge of the exact surface: the vertices of the triangles are assumed to be located on the exact surface, and the exact normal vectors are given.
However, this information is usually not available in reality, where we have only interpolation or approximations of surfaces, for instance, polyhedrons, splines or polynomial surfaces. How to deal with gradient recovery in such cases and prove its superconvergence is an open question in \cite{WeiChenHuang2010}.
On the other hand, the generalized $ZZ$-scheme  for surface elements gives the most competitive results in \cite{WeiChenHuang2010}, but its superconvergence, including several other methods, is proven with a condition that the mesh is $\mathcal{O}(h^2)$-symmetric. 
In the planar cases, this restrictive condition, however, is not necessary for the PPR method.

This triggers us to think of a generalization of the PPR method to manifolds setting.
A follow-up question would be what are  polynomials in the curved manifold domain.
Using the idea from the literature, e.g., \cite{DuJu2005}, one could consider polynomials defined locally on the tangent spaces of the manifolds.
Apparently, such a straightforward generalization based on tangent spaces will again fall into the awkward situation: the exact manifolds and their tangent spaces are unknown.

To overcome these difficulties, we go back to the original definition of manifolds. Every local patch of a manifold resembles a planar Euclidean domain, therefore a local parametrization for a patch of the manifold can always be established with respect to a parametric domain $\Omega$ and not necessarily to be a tangent space.
The idea is then to use polynomials to recover the unknown parametrization function of the discretized patch locally on the parametric domain $\Omega$, as well to fit the corresponding local data or finite element solutions on $\Omega$ iso-parametrically, from which we are able to recovery the gradient using the intrinsic definition (see formula \eqref{eq:local_gradient1} and \eqref{eq:local_gradient2}).
Our proposed method is called parametric polynomial preserving recovery (PPPR) which \emph{does not} rely on the $\mathcal{O}(h^2)-$ symmetric condition for the superconvergence, just like its genetic father PPR.
To this end, it is revealed that the idea of using parametric domain is particularly useful to \emph{address the issue of unavailable tangent spaces and vertices},
and thus it enables us to answer the open problem in \cite{WeiChenHuang2010}.
Even though we only prove the supperconvergence with exact vertices in the paper, as it can also be observed from our numerical examples, the superconvergence does hold for non-exact vertices of the triangulation.
To better demonstrate the ideas, a theoretical analysis of the problem with no exact vertices will be presented separately in a follow-up paper.
Another advantage of the PPPR method, which has been observed in our numerical examples, is that it is relatively \emph{high curvature stable} in comparing with the methods proposed in \cite{WeiChenHuang2010}.
This is verified by all of our numerical tests on the high curvature surfaces, but a quantitative analysis will be open in the paper.
Moreover, the original PPR method \cite{ZhangNaga2005} does not preserve the function values at the nodal points in its pre-recovery step.
In this paper, we take care of this issue, so that the PPPR can not only preserve \emph{parametric polynomial} but also preserve the \emph{surface sampling points} and the \emph{function values} at the given points simultaneously.
This property makes the given data invariant in using the PPPR method.

The rest of the paper is organized as follows:
Section \ref{sec:background} gives a preliminary account on relevant differential geometry concepts and an exemplary PDE problem.
Section \ref{sec:spaces} introduces discretized function spaces and collects some geometric notations frequently used in this paper.
Section \ref{sec:PPPR} presents the new algorithms especially the proposed PPPR method.
Section \ref{sec:analysis} exhibits the numerical analysis and some relevant properties of the PPPR method.
Section \ref{sec:estimator} shows the recovery-based {\it a posteriori} estimator by using the PPPR operator. 
Finally, in Section \ref{sec:numerics}, we present numerical results verifying the theoretical analysis and make comparisons with existing methods.
We postpone the proof of a basic lemma until Appendix \ref{appendix}.

\section{Background}
\label{sec:background}
We only show some basic geometric concepts which are relevant to our paper. For a more general overview on the topic of Riemannian geometry or differential geometry, one could refer to  \cite{doCarmo1992,Lee2013}.
In this paper, we shall consider $(\manifold,g)$ as an oriented, connected, $C^3$ smooth, regular and compact Riemannian manifold without boundary, where $g$ denotes the Riemann metric tensor. The idea we are going to work on should have no restriction for general $n$-dimensional manifolds, but we will focus on the case of two-dimensional ones, which are also called surfaces, in the later applications and numerical examples.

Our concerns are some quantities $u:\manifold \rightarrow \R$ which are scalar functions defined on manifolds.
First, let us recall the differentiation of a function $u$ in a manifold setting, which is called covariant derivatives in general.
It is defined as the directional derivatives of the function $u$ along an arbitrarily selected path $\gamma$ on the manifold
\[D_{\f v}u=\frac{du(\gamma(\sigma))}{d\sigma}|_{\sigma=0},\]
where $\f v=\gamma(\sigma)'|_{\sigma=0}$ is a tangential vector field.

The gradient then is an operator such that
\[(\nabla_g u(x),\f v(x))_g=D_{\f v}u,\; \text{ for all } \f v(x)\in T_x\manifold \text{ and all }  x\in \manifold,\]
where $T_x\manifold$ is the tangent space of $\manifold$ at $x$.
We can think  of the gradient as a tangent vector field on the manifold $\manifold$.
Using a local coordinate, the gradient has the form
\begin{equation}
\label{eq:local_gradient1}
\nabla_g u= \sum_{i,j} g^{ij}\partial_j u \partial_i,
\end{equation}
where $g^{ij}$ is the entry of the inverse of the metric tensor $g$, and $\partial_i$ denotes the tangential basis.
Fix $\f r:\Omega\rightarrow S \subset \manifold$ to be a local geometric mapping,
then we can rewrite \eqref{eq:local_gradient1} into a matrix form with this concrete local parametrization. That is
\begin{equation}
\label{eq:local_gradient2}
(\nabla_g u)\circ \f r=\nabla \bar{u}  (g\circ \f r)^{-1} \partial \f r.
\end{equation}
In \eqref{eq:local_gradient2}, $\bar{u}=u\circ \mathbf{r}$ is the pull back of function $u$ to the local planar parametric domain $\Omega$, $\nabla$ denotes the gradient operator on the planar domain $\Omega$, $\partial \f r$ is the Jacobian of $\f r$, and \[g\circ \f r =\partial \f r(\partial \f r)^\top.\]
\begin{remark}
\label{rem:surface_gradient}
$\f r$ is not specified here, and we will make it clear when it becomes necessary later.
We actually have a relation that 
\begin{equation}
\label{eq:jacob_inverse}
(\partial \f r)^\dag= (g\circ \f r)^{-1} \partial \f r,
\end{equation} 
where $(\partial \f r)^\dag$ denotes the Moore-Penrose generalized inverse of $\partial \f r$. See \cite[Appendix]{DonJueSchTak17} for a detailed explanation.
\end{remark}

We consider $\manifold$ a regular manifold in the paper. By regular we mean that the Jacobian of the parametrization $\partial \f r$ and their inverse $(\partial \f r)^\dag$ are functions with bounded norms in the space $W^{2,\infty}$ on $\Omega$.
 
Note that the parametrization map $\f r$ is non-unique. Typical ones can be constructed through function graphs on the parametric domain,
which will be used in our later algorithms.
We have the following lemma whose  proof is given in Appendix \ref{appendix}.
\begin{lemma}
\label{lem:invariant}
The gradient calculated by using \eqref{eq:local_gradient2} is invariant for different regular bijective parametrization functions $\f r$.
\end{lemma}

Let $\omega=dvol$ be the volume form on $\manifold$, and $\partial_j\;(j=1,\cdots,n)$ be the tangential bases, $T\manifold=\bigcup_{x\in \manifold}T_x\manifold$ be the tangent bundle which consists of all the tangent planes $T_x\manifold$ of $\manifold$. For every tangent vector field $\f v:\manifold\rightarrow T\manifold$, $\f v=v^i\partial_i$, we have a $(n-1)$ form defined by the interior product of $\f v$ and the volume form $\omega$ through the following way
\[i_{\f v}\omega=\sum_k \omega(\f v,\partial_{k_1},\cdots,\partial_{k_{n-1}}),\]
where $k_1,\cdots,k_{n-1}$ are $(n-1)$ indexes with $k$ taking out from $1,\cdots,n$.
The divergence of the vector field $\f v$  satisfies
\begin{equation}
\label{eq:divergence}
d(i_{\f v}\omega)=\text{div}_g(\f v) \omega,
\end{equation}
where $d$ denotes the exterior derivative.
Since both the left hand side and the right hand side of \eqref{eq:divergence} are $n$ forms, $\text{div}_g(\f v)$ is a scalar field.
Using the local coordinates, we can explicitly write the volume form as
\[\omega=\sqrt{\abs{\det g}}dx^1\wedge\cdots\wedge dx^n.\]
By equation \eqref{eq:divergence}, the divergence of the vector field $\f v $ can be computed by
\[\text{div}_g \f v=\frac{1}{\sqrt{\abs{\det g}}}\partial_i(v^i\sqrt{\abs{\det g}}).\]
It implies that the divergence operator is actually the dual of the gradient operator.
With the above preparation, we can now define the Laplace-Beltrami operator,
which is denoted by $\lapla_g$ in our paper, as the divergence of the gradient, that is
\begin{equation}
\label{eq:laplace_beltrami}
\lapla_g u=\text{div}_g(\nabla_g u)=\frac{1}{\sqrt{\abs{\det g}}}\partial_i(g^{ij}\sqrt{\abs{\det g}}\partial_j u).
\end{equation}
We would like to  mention that if the manifold $\manifold$ is a hypersurface, that is $\manifold \subset \R^{n+1}$ and it has co-dimension $1$, 
then the gradient and divergence of the function $u$ can be equally calculated through projecting the gradient and divergence of an extended function in ambient space $\R^{n+1}$ to the tangent spaces of $\manifold$ respectively.
That is
\[\nabla_g u=(\mathcal{P}_T \nabla_e) u_e \; \text{ and } \;\text{div}_g \f v=(\mathcal{P}_T \nabla_e)\cdot \f v_e,\]
where $u_e$ and $\f v_e$ are the extended scalar and vector fields defined in the ambient space of the hypersurface,
which satisfies $u_e(x)=u(x)$ and $\f v_e(x)=\f v(x)$ for all $x\in \manifold$.
Note that $\nabla_e$ is the gradient operator defined in the ambient Euclidean space $\R^{n+1}$,
$\mathcal{P}_T$ is the tangential projection operator
\[\mathcal{P}_T=\I-\f n\otimes \f n,\]
and $\f n$ is a unit normal vector field of $\manifold$.
Such type of definitions has been applied in many references e.g. \cite{WeiChenHuang2010} which consider problems in an ambient space setting. 

With the definition of covariant derivatives on manifolds, many function spaces on Euclidean domains can be studied analogously in the setting of  manifolds.
Sobolev spaces on manifolds \cite{Hebey1999} are one of the most investigated spaces, which provide a breeding ground to study PDEs.
We are interested in numerically approximating PDEs whose solutions are defined on $\manifold$.
Even though our methods are \emph{problem independent}, in this paper, the analysis will be mainly conducted for the Laplace-Beltrami operator \eqref{eq:laplace_beltrami} and its generated PDEs.
For the purpose of both analysis and applications, we consider the \emph{Laplace-Beltrami equation} as an exemplary problem \cite{Dziuk1988}:
For a given $f$ satisfying $\int_\manifold f \;dvol=0$,  find $u$ solves the equation
\begin{equation}
\label{eq:laplace}
-\lapla_g u= f \textit{   on   } \manifold,\; \text{  with  }  \int_\manifold u \;dvol=0,
\end{equation}
where $dvol$ denotes the manifold volume measure.

\section{Function Spaces on Discretized Manifolds}
\label{sec:spaces}
The discretization of a smooth manifold $\manifold$ has been widely studied in many settings, especially in terms of surfaces \cite{DziukElliott2013}.
A discretized surface, in most cases, is a piecewise polynomial surface.
One of the simplest cases is the polygonal approximation to a given smooth surface, especially with triangulations.
Finite element methods for triangulated meshes on surfaces have firstly been studied in \cite{Dziuk1988} by using the linear element.
In \cite{Demlow2009}, a generalization of \cite{Dziuk1988} to high order finite element methods is proposed based on triangulated surfaces. 
In order to have an optimal convergence rate, it is shown that the geometric approximation error and the function approximation error has to be compatible with each other.
In fact, the balance of the  geometric approximation error and the function approximation error is also the key point in the development of our recovery algorithm. 

For convenience, Table \ref{tab:Geometry} collects notations been frequently referred in the paper.
\begin{table}[!h]
\setlength{\tabcolsep}{3pt}
\caption{Notations }
\begin{center}
\begin{tabular}{ ll}
\hline
Notation & Remark    \\
\hline
$(\manifold,g)$  & a smooth, connected, oriented and close manifold with metric $g$        \\
\hline
$(\manifold_h ,g_h)$& a triangular approximation of $\manifold$ with metric $g_h$  \\
\hline
$\f n $& a unit normal vector field on $\manifold$    \\
\hline
$\nabla_g $& gradient operator with respect to the metric $g$ \\
\hline
$\lapla_g $& Laplace-Beltrami operator with respect to the metric $g$\\
\hline
$T_x$& a local domain on the tangent space at a position $x\in \manifold$    \\
\hline
$ (P_h)^{\pm 1} $ &  bijective maps between $\manifold_h$ and $\manifold$ \\
\hline
$ \mathcal{V}(\manifold) /\mathcal{V}_h(\manifold_h)$&  ansatz function spaces for functions on $\manifold$/ $\manifold_h $\\ 
\hline
$  (T_h)^{\pm 1}$ & operators between function spaces on $\manifold$ and on $\manifold_h$\\ 
\hline
$ h $ & the diameter of the triangulation mesh in $\manifold_h $  \\
\hline
$\Omega $ & a parametric domain for a patch on $\manifold/\manifold_h$   \\
\hline
$ \zeta $ & a position variable in the parameter domain $\Omega$  \\
\hline
$\f r /\f r_h $& a local parametrization map from $\Omega$ to a patch of $\manifold/ \manifold_h $  \\
\hline
$vol (\text{or } vol_h) $ & the volume (area) measure of $\manifold$ (or $\manifold_h$)   \\
\hline
$\norm{\cdot}_{k,p,\manifold} $ & $W^{k,p}$ norm of functions defined on $\manifold$  \\ 
\hline
$\abs{\cdot}_{k,p,\manifold}  $ & $W^{k,p}$ semi-norm of functions defined on $\manifold$   \\
\hline
$\norm{\cdot}_{k,\manifold}  $ & $H^k$ norm of functions defined on $\manifold$   \\
\hline
$I_h $ & the total number of the nodal points (vertices) of $\manifold_h$   \\
\hline
$J_h$ & the total number of the triangles on $\manifold_h$   \\
\hline
$\mathbb{P}_2(\Omega)$ & the $2^{nd}$ order polynomial space over a planar domain $\Omega$  \\
\hline
$a \circ b$  & function $a$ composed with function $b$  \\
\hline
$\alpha \lesssim \beta $ & denotes the inequality $\alpha\leq C \beta$ where $C$ is a constant  \\
\hline
$\mathcal{O}(\sigma)  $ & denotes the quantity satisfies: $\lim_{\sigma\to 0}\frac{\mathcal{O}(\sigma)}{\sigma} = C $ for $\sigma>0$ \\
\hline
\end{tabular}
\end{center}
\label{tab:Geometry}
\end{table}

Let $\manifold_h=\bigcup_{j\in J_h} \tau_{h,j}$ be a triangular mesh, and $h=\max_{j\in J_h} \mbox{diam}(\tau_{h,j})$ be the maximum diameter.
To better present our main idea, we mostly stick to the simplest case which is the linear finite elements on triangulated surfaces, thus the nodes consist of simply the vertices of $\manifold_h$, and we denote the set by $\mathcal{N}_h=\{x_i\}_{i\in I_h}$.

In the following, we define transform operators between the function spaces on $\manifold$ and on $\manifold_h$. Let $\mathcal{V}(\manifold)$ and $\mathcal{V}(\manifold_h)$ be some ansatz function spaces. Then we define
\begin{equation}
 \label{eq:transform}
 \begin{aligned}
  T_h : \mathcal{V}(\manifold) &\to  \mathcal{V}(\manifold_h); \\
        v &\mapsto  v \circ P_h,
\end{aligned}
\end{equation}
and its inverse
\begin{equation}
 \label{eq:itransform}
 \begin{aligned}
  (T_h)^{-1}:  \mathcal{V}(\manifold_h) &\to  \mathcal{V}(\manifold); \\
      v_h &\mapsto   v_h \circ P_h^{-1},
\end{aligned}
\end{equation}
where $P_h$ is a continuous and bijective projection map from every element in $\set{\tau_{h,j}}_{j\in J_h}$ to every element in $\set{\tau_j}_{j\in J_h}$.

We will use the following definition to characterize the approximation quantity of $\manifold_h$ to $\manifold$.
For the purpose of later analysis, we assume that both $T_{x_i}$ and $\Omega_i$ are compact domain corresponding to selected compact patches on $\manifold_h$ or $\manifold$.
\begin{definition}
\label{def:h2surface_appr}
Let $\manifold_h=\bigcup_{j\in J_h} \tau_{h,j}$ be a triangular approximation of $\manifold$.
Let $\mathcal{K}_i\subset \manifold_h$ be the triangle patches associated with the vertex $x_i\in \mathcal{N}_h$.  Then there is a curved patch  $\manifold_i\subset \manifold$.
Let $\mathcal{K}_i$ and $\manifold_i$ be parametrizable by a common domain $\Omega_i$ with $\mathbf{r}_{h,i}$ and $\mathbf{r}_i$ be their parametrization functions respectively.
We call $\manifold_h$ is a \emph{regular approximation} of $\manifold$ if
\begin{equation}
\label{eq:surface_converge}
\lim_{h\rightarrow 0}\norm{\mathbf{r}_{h,i}-\mathbf{r}_i}_{\infty;\Omega_i} =0\;  \text{ for all } i\in I_h.
\end{equation}
for a fixed number $k\in \N$, and both $\abs{\partial \mathbf{r}_{h,i}} $ and its inverse $\abs{(\partial \mathbf{r}_{h,i})^\dag} $ are uniformly bounded on $\Omega_i$ (edges are ignored) for all $i\in I_h$.
\end{definition}

Based on our assumptions on $\manifold$, it indicates that  $\f r_i\in W^{3,\infty}(\Omega_i)$ for all $i\in I_h$.
In particular, since we consider the vertices of $\manifold_h$ located on $\manifold$, therefore, $\f r_{h,i}$ is a linear interpolation of $\f r_i$.
If $\manifold_h$ is a regular approximation of $\manifold$, then it converges to $\manifold$ as $h\rightarrow 0$.
In the following, we introduce conditions on the triangle meshes which are common conditions to guarantee  the  supercloseness (cf. \cite[Definition 2.4]{BankXu2003}, \cite[Definition 1.2]{NagaZhang2004} or \cite[Definition 3.2]{WeiChenHuang2010}).

\begin{figure}[htbp]
\begin{center}
   \includegraphics[width=0.3\textwidth]{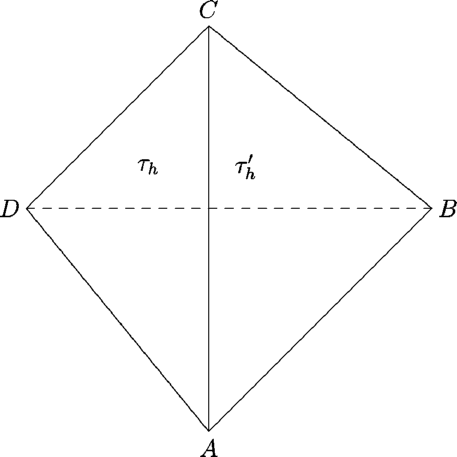}
\caption{Illustration of two adjacent triangles.}
\label{fig:parallel}
\end{center}
\end{figure}

\begin{definition}
\label{def:parallel}
Suppose $\tau_h$ and $\tau_h'$ are two adjacent triangles  in $\mathcal{T}_h$, as illustrated in Figure \ref{fig:parallel}. 
They are said to form  an $\mathcal{O}(h^2)$ parallelogram if 
\begin{equation*}
|\overrightarrow{AB}-\overrightarrow{CD}|  = \mathcal{O}(h^2), 
\quad\text{ and }\quad
|\overrightarrow{BC}-\overrightarrow{DA}|  = \mathcal{O}(h^2).
\end{equation*}

\end{definition}

\begin{definition}
\label{def:2sigma_irregular}
A triangulation mesh $\mathcal{T}_h$ is said to satisfy the \emph{$\mathcal{O}(h^{2\sigma})$ irregular condition} if
there exist a partition $\mathcal{T}_{h,1} \bigcup \mathcal{T}_{h,2}$ of $\mathcal{T}_h$ and a positive constant $\sigma$ such that every two adjacent triangles in $\mathcal{T}_{h,1} $ form an $\mathcal{O}(h^{2})$ parallelogram and 
\[\sum_{\tau_h\subset \mathcal{T}_{h,2}} \abs{\tau_h} =\mathcal{O}(h^{2\sigma}).\]
\end{definition}

Before going further, we make a general assumption for a non-adaptive triangulation $\manifold_h$.
\begin{assumption}
\label{ass:irregular}
Let $\manifold_h$ be a triangulation of $\manifold$ with all of the nodes located on $\manifold$. We assume it to be quasi-uniform and shape regular, and be a regular approximation of $\manifold$. Moreover, it satisfies the $\mathcal{O}(h^{2\sigma})$ irregular condition.
\end{assumption}

We have the following lemma for the transform operators.
\begin{lemma}
\label{lem:transform}
Let $\mathcal{V}(\manifold) \hookrightarrow W^{k,p}(\manifold)$ for a fixed $k\in \N$ and $ p\geq 1$, and $\manifold_h$ satisfies Assumption \ref{ass:irregular}. Then the predefined operators $(T_h)^{\pm 1}$ are uniformly bounded between the spaces $W^{k,p}(\manifold)$ and $W^{k,p}(\manifold_h)$, that is there exists positive constants $c\leq C$, and
\[ \quad c\norm{T_h v}_{k,p,\manifold_h}  \leq \norm{v}_{k,p,\manifold}\leq C  \norm{T_h v}_{k,p,\manifold_h} \text{ for all } v\in \mathcal{V}(\manifold).\]
\end{lemma}
\begin{proof}
Denote $\check{v}_h:=T_h v$, and let $\tau_j\subset \manifold$ be the curved triangle corresponding to $\tau_{h,j}\subset \manifold_h$. If $p=\infty$, every function $v$ and its derivatives are uniformly bounded on $\manifold$, as well for function $\check{v}_h$ and its derivatives over $\manifold_h$. Then we can always find constants $c^1_{h} $ and $C^1_{h}$ satisfying 
\[c^1_{h} \norm{\check{v}_h}_{k,\infty,\manifold_h} \leq \norm{v}_{k,\infty,\manifold} \leq C^1_{h}\norm{\check{v}_h}_{k,\infty,\manifold_h}. \]
If $1\leq p <\infty$, using the results in \cite[page 811]{Demlow2009}, there exists positive and bounded constants $c_{h,j}$ and $C_{h,j}$ for each pair of triangle faces $\tau_{h,j}$ and $\tau_j$ such that 
\[c^2_{h,j} \norm{\check{v}_h}_{k,p,\tau_{h,j}}^p\leq \norm{v}_{k,p,\tau_j}^p \leq C^2_{h,j}\norm{\check{v}_h}_{k,p,\tau_{h,j}}^p.\]
For both the two cases, due to the regular approximation condition in Assumption \ref{ass:irregular}, we have $c^1_h\rightarrow 1$, $C^1_h\rightarrow 1$ when $h\rightarrow 0$ as well as $c^2_{h,j} \rightarrow 1$ and $C^2_{h,j}\rightarrow 1$ when $h\rightarrow 0$ for all $j\in J_h$ .
Thus, both $\set{c^a_{h,j}}$ and $\set{C^a_{h,j}}$ are uniformly bounded sequences with respect to the mesh size $h$ and also the index $j$ for $a=1,2$. Denote $c:=\min_{a,h,j} \set{c^a_{h,j}}$ and $C:=\max_{a,h,j} \set{C^a_{h,j}}$.
Since $\norm{v}_{k,p,\manifold}^p=\sum_{j\in J_h} \norm{v}_{k,p,\tau_j}^p$ for $p\in [1,\infty)$,  we have  the estimates
\[c \norm{\check{v}_h}_{k,p,\manifold_h}^p \leq \norm{v}_{k,p,\manifold}^p \leq C \norm{\check{v}_h}_{k,p,\manifold_h}^p,\] 
and the same hold for $p=\infty$. This gives the conclusion.
\end{proof}


\section{Parametric Polynomial Preserving Recovery on Manifolds}
\label{sec:PPPR}
Our developments are based on the  PPR method proposed in \cite{ZhangNaga2005} for planar problems.
It is a robust and high accuracy approach for recovering gradient on mildly unstructured  meshes. This idea has been used to develop a Hessian recovery technique in a recent paper \cite{GuoZhangZhao2016}.
In this paper, we show the possibility of generalizing the idea to problems on manifolds. To simplify the presentation, we shall restrict ourselves to the case of two-dimensional manifolds here and after.

We will focus on the case where the data is a linear finite element solution on $\manifold_h$.  Therefore,  $\mathcal{V}_h(\manifold_h)$ is restrict to linear finite element spaces in what follows.
At each node $x_i$,  let $h_i$ be the length of the longest edge attached to $x_i$.
For any natural number $k$, let $B_{kh_i}(x_i)$ be the set of vertices in a discrete geodesic ball centered at $x_i$ with discrete geodesic radius $k\times h_i$, i.e., 
\begin{equation*}
B_{kh_i}(x_i) = \{x\in \mathcal{N}_h: |x-x_i|\le k\times h_i\}.
\end{equation*}
Then we define $B(x_i) = B_{k_ih_i}(x_i)$ with $k_i$ being the smallest integer such that $B(x_i)$ satisfies the rank 
condition (see \cite{ZhangNaga2005}) in the following sense:
\begin{definition}\label{def:rank}
 A selected vertices set $B(x_i)$ is said to satisfy the rank condition of the PPR or PPPR if it admits
 a unique least-squares fitted polynomial $p_i$ in \eqref{equ:alg1ls} or $s_i$ and $p_i$ in \eqref{equ:alg22s} and \eqref{equ:alg2ls} respectively.
\end{definition}

For a discretized manifold, the main difficulty is that the vertices in $B(x_i)$ are in general not located on the same plane. Another challenge is that there is no trivial definition of polynomials in a manifold setting. Some idea appeared in the literature is to use the domain in tangent space $T_{x_i}$ at every vertex $x_i$ as a local parameter domain, and project the neighboured vertices of $x_i$ onto this common planar plane, then define polynomials locally by the coordinates of the tangent space. This idea has been applied in \cite{DuJu2005} and also in \cite{WeiChenHuang2010} to generalize the $ZZ$ method and several other methods.
However, the exact manifold $\manifold$ is usually not given in real problems. Therefore, the tangent spaces $(T_{x_i})_{i\in I_h}$ of  $\manifold$ are blind to users, which makes the idea not much feasible in practice. This problem has been proposed as an open question in \cite{WeiChenHuang2010}. 

In our algorithms, we consider an alternative way for the polynomial reconstruction instead of the one which is initially proposed in \cite{ZhangNaga2005} for the  planar PPR. 
The method in \cite{ZhangNaga2005} assumes that a second order polynomial has a form 
\[p(y)=a_0+a_1y_1+a_2y_2+a_3y_1^2+a_4y_1y_2+a_5y_2^2 , \text{ for } \; y=(y_1,y_2)\in \Omega_i,\] 
and solves the linear system $A \f a=\f b$ for $\f a=(a_0,a_1,\cdots,a_5)^T$, where
\begin{equation}
\label{eq:matrix1}
A=
\left(
\begin{matrix}
1\;& \zeta_{i_1,1} \;& \zeta_{i_1,2}\;& \zeta_{i_1,1}^2\;& \zeta_{i_1,1}\zeta_{i_1,2}\; &\zeta_{i_1,2}^2\\
1\;& \zeta_{i_2,1} \;& \zeta_{i_2,2}\;& \zeta_{i_2,1}^2\;& \zeta_{i_2,1}\zeta_{i_2,2}\; &\zeta_{i_2,2}^2\\
\vdots &\cdots\\
1\;& \zeta_{i_{\abs{I_i}},1} \;& \zeta_{i_{\abs{I_i}},2}\;& \zeta_{i_{\abs{I_i}},1}^2\;& \zeta_{i_{\abs{I_i}},1}\zeta_{i_{\abs{I_i}},2}\; &\zeta_{i_{\abs{I_i}},2}^2\\
\end{matrix}
\right) 
\text{ and } 
\f b=
\left(
\begin{matrix}
u_{h,i_1}\\
u_{h,i_2}\\
\vdots\\
u_{h,i_{\abs{I_i}}}\\
\end{matrix}
\right).
\end{equation}
The solution of the least squares approximation in the algorithms is given by
\[\f a = (A^T A)^{-1}  A^T \f b, \]
which tells that $\partial_1 p(0)=a_1$ and $\partial_2 p(0)=a_2$.

Our observation is that there is one extra freedom can be removed in the reconstruction of the polynomials.
Since the polynomial recovery procedure cannot improve the accuracy of the solution itself, it is unnecessary to adopt the solution in gradient recovery. We can fix this problem by using the following polynomial equation locally
\[\tilde{p}(y)=u_{h,i_1}+\tilde{a}_1 y_1+\tilde{a}_2 y_2+\tilde{a}_3 y_1^2+\tilde{a}_4y_1 y_2+\tilde{a}_5 y_2^2, \text{ for } \; y=(y_1,y_2)\in \Omega_i\]
where $u_{h,i_1}$ is the finite element solution at the vertex $x_i$.
Let $\zeta_{i_1}=(\zeta_{i_1,1},\zeta_{i_1,2})$ be the origin $0$ of the plane $\Omega_i$, then the matrix and the vector in \eqref{eq:matrix1} can be simplified to
\begin{equation}
\label{eq:matrix2}
\tilde{A}=
\left(
\begin{matrix}
\zeta_{i_2,1} \;& \zeta_{i_2,2}\;& \zeta_{i_2,1}^2\;& \zeta_{i_2,1}\zeta_{i_2,2}\; &\zeta_{i_2,2}^2\\
\vdots&\cdots\\
\zeta_{i_{\abs{I_i}},1} \;& \zeta_{i_{\abs{I_i}},2}\;& \zeta_{i_{\abs{I_i}},1}^2\;& \zeta_{i_{\abs{I_i}},1}\zeta_{i_{\abs{I_i}},2}\; &\zeta_{i_{\abs{I_i}},2}^2\\
\end{matrix}
\right)
\text{ and } 
\tilde{\f b}=
\left(
\begin{matrix}
u_{h,i_2}-u_{h,i_1}\\
\vdots\\
u_{h,i_{\abs{I_i}}}-u_{h,i_1}\\
\end{matrix}
\right).
\end{equation}
Solving the problem in the least squares sense
\[\tilde{\f a} =(\tilde{A}^T\tilde{A})^{-1} \tilde{A}^T \tilde{\f b}, \]
then we have $\partial_1 \tilde{p}(0)=\tilde{a}_1$ and $\partial_2 \tilde{p}(0)=\tilde{a}_2$.

Using \eqref{eq:matrix2} instead of \eqref{eq:matrix1}, it preserves the original function values at the recovered nodal points. This idea can be applied to construct the polynomial functions in the $2^{nd}$ step of Algorithm \ref{alg:1}, and also in both the $2^{nd}$ and the $3^{rd}$ steps of Algorithm \ref{alg:2}, which are going  to be introduced next.

As a starting point, we first provide a direct generalization of the PPR method based on given tangent spaces of the exact manifold $\manifold$. In this case, the algorithm is pretty much  the same as the planar one. We sketch it in Algorithm \ref{alg:1} and still name it  as the PPR method. We describe the PPPR method in Algorithm \ref{alg:2}. In both algorithms, $I_{h,i}$ denotes the set of the indexes of the selected vertices in $B(x_i)$ which satisfies the rank condition.   

\begin{algorithm}
\caption{PPR Method (\emph{with exact information of normal vectors})}
\label{alg:1}
\begin{minipage}{\textwidth}
Let the discretized triangulation $\manifold_h$ and the data (FEM solutions) $(u_{h,i})_{i\in I_h}$ be given.
Also, we have the the normal vector $(\f n_i)_{i\in I_h}$ of $\manifold$ at each vertex $x_i$.
Then repeat steps $(1)-(3)$ for all $i\in I_h$.
\begin{itemize}
\item[(1)] 
For every $x_i$, select $B(x_i)\in \manifold_h$ including sufficient vertices, and shift $x_i$ to be the origin of $T_{x_i}$, and choose an orthonormal basis $(\f t^1_i, \f t^2_i)$ of $T_{x^*_i}$, then project the vertices $x^*_j\in B(x_i)$ to $T_{x_i}$ whose new coordinates read as $\zeta_{i_j}$. 
\item[(2)] Find a polynomial $p_i$ over $T_{x_i}$ by solving the least squares problem
\begin{equation}\label{equ:alg1ls}
p_i=\arg\min_{p} \sum_{j\in I_{h,i}} \abs{p(\zeta_{i_j})-u_{h,j}}^2\; \text{ for } p\in \mathbb{P}_2(T_{x_i}).
\end{equation}

\item[(3)] Calculate the partial derivatives of the approximated polynomial functions to have the recovered gradient at each vertex $x_i$
\begin{equation}
\label{eq:gradient_ppr}
G^*_{h} u_h(x_i)=\partial_1 p_i(0)\f t^1_i+\partial_2 p_i(0)\f t^2_i.
\end{equation}
\end{itemize}
For the recovery of  gradient $G^*_{h} u_h(x)$ when $x$ is not a vertex of triangles, use linear finite element basis to interpolate the values $\set{G^*_{h} u_h(x_i)}_{i\in I_h}$ at vertices of each triangle.
\end{minipage}
\end{algorithm}

A straightforward remedy for missing exact normal fields is to find a way to approximate normal vectors at every vertex $x_i$. This can be done, for instance, by  taking the simple average or weighted average of the normal vectors of each faces adjunct to $x_i$.
However, with such kinds of  approximations, the recovery errors are very likely to be dominated by the errors of the approximation of the normal vector fields (see the numerical results in Section \ref{sec:numerics}).
Therefore,  a better estimation of the normal vectors is  necessary in order to have higher recovery accuracy.

\newpage
The PPPR method (Algorithm \ref{alg:2}) requires no information of the tangent spaces of $\manifold$.
The idea is to use the intrinsic formulation \eqref{eq:local_gradient2}, where we can calculate the gradient from an arbitrary local parametrization.
The local parametrization function can be constructed, in principal, with respect to arbitrary Euclidean domain $\Omega_i$ but not restrict to tangent space $T_{x_i}$ in Algorithm \ref{alg:1}.
Lemma \ref{lem:invariant} indicates that for every fixed $x_i$, taking arbitrary $\Omega_i$, the gradient operator is analytically invariant. 
The crucial point in practice is that, numerically, the shape of the triangles must not be destroyed after projecting them to the domain $\Omega_i$, and also for a superconvergence purpose, the $\mathcal{O}(h^{2\sigma})$ irregular condition should be properly preserved for the projected triangular mesh on $\Omega_i$.
Thus, we still have to find a good way for this projection.
Our suggestion is, at each vertex, to use the simple average or weighted average of the surrounding normal vectors which can help us to locate and orient a suitable parameter domain $\Omega_i$.
This has been adopted in our numerical examples.
\begin{algorithm}
\caption{PPPR Method (\emph{with no exact normal vectors})}
\label{alg:2}
Let the discretized triangular surface $\manifold_h$ and the data (FEM solutions) $(u_{h,i})_{i\in I_h}$ be given.
Then repeat steps $(1)-(4)$ for all $i\in I_h$.
\begin{itemize}
\item[(1)] For every $x_i$, select $B(x_i)\in \manifold_h$ with sufficient vertices, using simple (weighted) average of the out normal vectors of every triangles with vertices in $B(x_i)$, and normalizing the averaged vector to be $\phi^3_i$, and then constructing a local parameter domain $\Omega_i$ orthogonal to $\phi^3_i$. Shift $x_i$ to be the origin of $\Omega_i$, and choose $(\phi^1_i,\phi^2_i)$ the orthonormal basis of $\Omega_i$, then project all selected vertices $x_j\in B(x_i)$ into the parameter domain $\Omega_i$, and record the new coordinates as $\zeta_{i_j}$. 

\item[(2)] Reconstruct a $2^{nd}$ order polynomial surface $S_i$ over $\Omega_i$ to approximate the local surface.
Typically, it can be approximated locally as a function graph parametrized by $\Omega_i$. That is $S_i=\tilde{\f r}_{h,i}(\Omega_i)=\bigcup_{\zeta\in \Omega_i}(\zeta,s_i(\zeta))$, where $s_i$ solves
\begin{equation}\label{equ:alg22s}
s_i=\arg\min_{s} \sum_{j\in I_{h,i}} \abs{s(\zeta_{i_j})-\langle x_j,\phi^3_i \rangle }^2\; \text{ for } s\in \mathbb{P}_2(\Omega_i).
\end{equation}
\item[(3)] Find a $2^{nd}$ order polynomial $p_i$ over the domain $\Omega_i$ by optimizing
\begin{equation}\label{equ:alg2ls}
p_i=\arg\min_{p} \sum_{j\in I_{h,i}} \abs{p(\zeta_{i_j})-u_{h,j}}^2\; \text{ for } p\in \mathbb{P}_2(\Omega_i).
\end{equation}

\item[(4)] Calculate the partial derivatives of both the polynomial approximated surface function in Step (2) and the approximated polynomial function of FEM solution in Step (3). Then approximate the gradient as \eqref{eq:gradient_pppr} using the local coordinates: 
\begin{equation}
\label{eq:gradient_pppr}
G_h u_h(x_i) =
\left(
\begin{matrix}
\partial_1 p_i(0),
\partial_2 p_i(0)
\end{matrix} 
\right)
J^\dag(s_i)  
\left(
\begin{matrix}
\phi^1_i \; \phi^2_i\; \phi^3_i
\end{matrix}
\right)^\top,
\end{equation}
where $J^\dag(s_i)=(J^\top(s_i) J(s_i))^{-1}J^\top(s_i)$ and $J^\top(s_i)=\left(
\begin{matrix}
1 & 0& \partial_1 s_i(0)\\
0 & 1& \partial_2 s_i(0)
\end{matrix}
\right)$. The equation \eqref{eq:gradient_pppr} is derived from \eqref{eq:jacob_inverse} in the remark \ref{rem:surface_gradient} for calculating \eqref{eq:local_gradient2}.
To multiply with the orthonormal basis $\set{\phi^1_i, \phi^2_i, \phi^3_i}$ is because we have to unify the coordinates from local ones to a global one.
\end{itemize}
For the recovery of the gradient $G_{h} u_h(x)$ when $x$ is not a vertex of triangles, use linear finite element basis to interpolate the values $\set{G_{h} u_h(x_i)}_{i\in I_h}$ at vertices of each triangle.
\end{algorithm}
Note that for $\Omega_i=T_{x_i}$ for all $i\in I_h$, and we shift $x_i$ to be the origin of $T_{x_i}$, 
and let $\phi^1_i=\f t^1_i$, $\phi^2_i=\f t^2_i$, $\phi^3_i=\f n_i$,
if $\partial_1 s_i(0)=\partial_2 s_i(0)\equiv 0$ for all $i\in I_h$, then the recovered gradient in \eqref{eq:gradient_pppr} is equal to the one recovered in \eqref{eq:gradient_ppr}.

Let $\bar{G}_h$ be the PPR operator introduced in  \cite{NagaZhang2004} for planar problems.
The recovered gradient values at each vertex $x_i$ by PPPR operator $G_h$ can be represented by $\bar{G}_h$ in the following sense:
\begin{equation}
\label{eq:PPR_relations}
G_h u_h(x_i) =\bar{G}_h \bar{u}_h(\zeta_i) (\bar{G}_h \f r_{h,i}(\zeta_i))^\dag, \quad \zeta_i \in \Omega_i \text{ is the projection from } x_i.
\end{equation}

Our numerical results will show that choosing the approximations of normal vectors by either simple average or weighted average has very little influence on the recovery accuracy of the gradient by Algorithm \ref{alg:2}.
This is different to the case in Algorithm \ref{alg:1} where the recovery accuracy highly relies on the error of the approximated normal vectors. The relation \eqref{eq:PPR_relations} indicates that the analysis of the PPR which has been developed for planar problems can be applied to Algorithm \ref{alg:2} to some extend.
Moreover, the idea of approximating \eqref{eq:local_gradient2} by generalizing $ZZ$ scheme seems feasible.
One could similarly reconstruct the two levels gradient recovery of the surfaces parametrization function $\mathbf{r}$ and  the function $\bar{u}$ iso-parametrically. That is to replace the recovery operator $\bar{G}_h$ in \eqref{eq:PPR_relations} by using planar $ZZ$ recovery.
However, in order to achieve the superconvergence property, this generalization can still not escape the constraint that the meshes should be $\mathcal{O}(h^2)-$ symmetric.

\begin{remark}
 \label{rem:highcurvature}
An experimental observation will be reported later that the PPPR is able to give the most competitive results for the recovery of the gradient when the approximated surface is featured with some high curvature.
Our argument is that, in the planar case, the PPR is the most robust method for unstructured meshes  compared to the other methods, especially, it does not require the $\mathcal{O}(h^2)$ symmetric condition.
For a surface with complicated curvature, a well-structured triangulation after projecting to the parametric domains or tangent spaces, it is possible that the good structure is not preserved, for instance,   the symmetric condition.
The PPPR method is, in fact, using the PPR to reconstruct both the tangent vectors of the surface and the gradient of the solutions in local parametric domains, which is more stable than the other methods for those mildly structured meshes projected from the high curvature areas.
All of our numerical tests on high curvature surfaces support this hypothesis, as the one shown in Numerical Example $2$ in Section \ref{sec:numerics}.
However, more efforts are needed in order to have a quantitative analysis of this property.
\end{remark}

\section{Superconvergence Analysis}
\label{sec:analysis}
We prove the superconvergence property of the proposed algorithms in previous  section.
Although our algorithms are problem independent,  to make the discussion simple, we will take the equation \eqref{eq:laplace} as our model problem, and focus on its  approximation using the linear finite element method on triangulated surfaces. 
The variational formulation of problem \eqref{eq:laplace} is given as follows: Find $u\in H^1(\manifold)$ such that
\begin{equation}
\label{eq:var_laplace}
\int_{\manifold} \nabla_g u \cdot \nabla_g v \; dvol=\int_{\manifold} f v \; dvol,  \text{  for all } \;v\in H^1(\manifold).
\end{equation}
The regularity of the solutions has been proved in \cite[Chapter 4]{Aubin1982}. In the finite element methods, the surface $\manifold$ is approximated by the triangulation $\manifold_h$ which satisfy Assumption \ref{ass:irregular},
and the finite element space $\mathcal{V}_h(\manifold_h)$ is   the piecewise linear function spaces  defined over $\manifold_h$. 
The finite element solution is to find $u_h \in \mathcal{V}_h(\manifold_h)$ such that 
\begin{equation}
\label{eq:d_var_laplace}
\int_{\manifold_h} \nabla_{g_h} u_h \cdot \nabla_{g_h}  v_h \; dvol_h=\int_{\manifold_h} f_h v_h \; dvol_h, \text{  for all } \;v_h\in \mathcal{V}_h(\manifold_h).
\end{equation}
Since here quite a few local geometric notations involved, we summarize them in Table \ref{tab:local_geometry}.
\begin{table}[!h]
\setlength{\tabcolsep}{3pt}
\caption{Notations on local geometry}
\begin{center}
\begin{tabular}{ ll}
\hline
Notation & Remark    \\
\hline
$\Omega_i$ & a local parametric domain for patches around vertex $x_i$   \\
\hline
$\mathcal{K}_i$ & local triangle patches of $\manifold_h$ around vertex $x_i$   \\
\hline
$\tau_{h,j}$  &  $\tau_{h,j} \subset \manifold_h$ is the $j^{th}$ triangle face \\
\hline
$\tau_{j}$  & $\tau_{j}\subset \manifold $ is the $j^{th}$ curved triangle w.r.t. $\tau_{h,j}\subset \manifold_h$   \\
\hline
$\tilde{\f r}_{\tau_{h,j}}$ & geometric mapping from $\tau_{h,j}$  to $\tau_{j}$    \\
\hline
$\f r_{h,i}$ & geometric mapping from patches in $\Omega_i$ to $\mathcal{K}_i$    \\
\hline
$\f r_{i}$ & geometric mapping from patches in $\Omega_i$ to $\manifold_i$    \\
\hline
\end{tabular}
\end{center}
\label{tab:local_geometry}
\end{table}

We prepare the proof of the superconvergence of the recovered gradient by firstly establishing the boundedness of the proposed gradient recovery operators. 
\begin{lemma}
\label{lem:boundedness}
Choose an arbitrary but fixed vertex $x_i$, and let $\tau_{h,j}$ be one of the triangles connected to $x_i$, $\mathcal{K}_i$ be the selected triangle patches, and $\tau_{h,j}\subset \mathcal{K}_i \subset \manifold_h$.
Then $G_h$ is a bounded linear operator in the sense that
\begin{equation}
\label{eq:boundedness}
\norm{G_h v_h}_{L^2(\tau_{h,j})}\lesssim \norm{\nabla_{g_h} v_h}_{L^2(\mathcal{K}_i)} ,\quad \text{ for all } v_h\in \mathcal{V}_h(\manifold_h).
\end{equation}
\end{lemma}
\begin{proof}
Let us denote $\bar{v}_h=v_h\circ \f r_{h,i}$, and recall \eqref{eq:local_gradient2}. Then we have that on every parametric domain $\Omega_i$:
\begin{equation}
\label{eq:invergradient}
(\nabla_{g_h} v_h)\circ \f r_{h,i}= \nabla \bar{v}_h (\partial\mathbf{r}_{h,i})^\dag \; \Longleftrightarrow  \;\nabla  \bar{v}_h =  (\nabla_{g_h} v_h)\circ\mathbf{r}_{h,i}  \partial\mathbf{r}_{h,i},
\end{equation}
where $\partial\mathbf{r}_{h,i}$ and $(\partial\mathbf{r}_{h,i})^\dag$ are piecewise constant functions.
We take into account the assumptions that $\manifold$ is regular and $C^3$ smooth, and $\manifold_h$ is a regular approximation as specified in Definition \ref{def:h2surface_appr}. 
Then there exist positive constants $c_r$ and $C_r$, such that
\begin{equation}
\label{eq:bounds}
 \frac{n}{C_r}\leq \abs{\partial\mathbf{r}_{h,i}}\leq \frac{n}{c_r},\;\text{  and   }\; c_r\leq \abs{(\partial\mathbf{r}_{h,i})^\dag} \leq C_r\;\text{ for all $h$ and $i$},
\end{equation}
where $n$ is the dimension number of $\manifold$. Correspondingly there are similar bounds for $\abs{(\partial \f r_i)^\dag}$ and we denote them with $c^*_r$ and $C^*_r$.
Because of the interpolation nature, and the fact that every $\tau_{h,j}$ is uniformly bounded, we get
\begin{equation}\label{eq:point_wise}
 \abs{G_h v_h}(x) \leq \sum_{i\in V_j}{\abs{G_h v_h(x_i)}} \leq C \max_{i\in V_j}{\abs{G_h v_h(x_i)}} \; \text{ for all } x\in \tau_{h,j};
\end{equation}
where $V_j$ denotes the index set of the vertices on $\tau_{h,j}$.
Using the boundedness result of the planar PPR recovery operator \cite{NagaZhang2004},  we have for every $\tau_{h,j}\subset\mathcal{K}_i$:
\[\norm{\bar{G}_h \bar{v}_h }_{L^\infty(\tau_{h,j})} \leq C\norm{\nabla  \bar{v}_h}_{L^\infty(\Omega_i)} .\]
Now we are going to show that $\norm{(\bar{G}_h \f r_{h,i})^\dag }_{L^\infty(\tau_{h,j})} $ is uniformly bounded for all $i\in I_h$ and $j\in J_h$. 
By noticing that $\f r_{h,i}$ is a linear interpolation of $\f r_i$ on $\Omega_i$,  the polynomial preserving property of the planar PPR operator implies
\begin{equation}
\label{eq:re_geo_error_bounds}
\abs{\bar{G}_h \f r_{h,i} -\partial \f r_i}=\norm{\f r_i}_{3,\infty}\mathcal{O}(h^2) .
\end{equation}
Taking into account that $\partial \f r_i$ is uniformly bounded from below and above and $\f r_i$ belongs to $W^{3,\infty}(\Omega_i)$, 
we can deduce from  \eqref{eq:re_geo_error_bounds} that there exists   a constant $c$ for sufficiently small $h^2\leq h^2_0$ (some fixed $h_0\in \R^+$) such that 
\begin{equation}\label{eq:lowerbound}
	\abs{\bar{G}_h \f r_{h,i}} \geq \frac{n}{C^*_r}- c h_0^2 \; \Rightarrow   \; \norm{(\bar{G}_h \f r_{h,i})^\dag }_{L^\infty(\tau_{h,j_i})}\leq \hat{C} ,
\end{equation} 
where $\hat{C} :=\frac{n}{\frac{n}{C^*_r}- c h_0^2 }$. For   every $\tau_{h,j}\subset \mathcal{K}_i$, we conclude that 
\[ \max_{i\in V_j}{\abs{G_h v_h(x_i)}} \leq \norm{\bar{G}_h \bar{v}_h }_{L^\infty(\tau_{h,j})}  \norm{(\bar{G}_h \f r_{h,i})^\dag }_{L^\infty(\tau_{h,j})}  \leq \hat{C}\norm{\nabla  \bar{v}_h}_{L^\infty(\Omega_i)};  \]
which together with \eqref{eq:point_wise} give
\[  \norm{G_h v_h}_{L^\infty(\tau_{h,j})}\leq C\hat{C}\norm{\nabla  \bar{v}_h}_{L^\infty(\Omega_i)}.\]
Using the formula on the right side of \eqref{eq:invergradient}, 
and the bounds on $\abs{\partial\mathbf{r}_{h,i}}$ in \eqref{eq:bounds},
 we get the boundedness result for $G_h$
\[\norm{G_h v_h}_{L^\infty(\tau_{h,j})}  \leq \frac{nC\hat{C}}{c_r} \norm{ \nabla_{g_h} v_h}_{L^\infty(\mathcal{K}_i)} .\]
All the constants here are independent of $h$ for all $h\leq h_0$.
Since $v_h$ is a piecewise linear polynomial on $\manifold_h$,  the inverse estimate implies
\[\norm{ \nabla_{g_h} v_h}_{L^\infty(\mathcal{K}_i)} \leq \frac{C_{in}}{\sqrt{ \abs{\mathcal{K}_i}}} \norm{ \nabla_{g_h} v_h}_{L^2(\mathcal{K}_i)}, \]
for some constant $C_{in}$ independent of $h$. Here $ \abs{\mathcal{K}_i}$ denotes the area of $\mathcal{K}_i$. 
Finally we have
\begin{eqnarray*}
\norm{G_h v_h}_{L^2(\tau_{h,j})}
 &\leq & \sqrt{\abs{\tau_{h,j}}}\norm{G_h v_h}_{L^\infty(\tau_{h,j})} \leq \sqrt{\abs{\tau_{h,j}}}\frac{nC \hat{C}}{c_r} \norm{ \nabla_{g_h} v_h}_{L^\infty(\mathcal{K}_i)} \\
 &\leq &  \frac{nC \hat{C}C_{in} \sqrt{\abs{\tau_{h,j}}}}{c_r\sqrt{ \abs{\mathcal{K}_i}}}\norm{ \nabla_{g_h} v_h}_{L^2(\mathcal{K}_i)} .
\end{eqnarray*}
The fact  $\frac{\sqrt{\abs{\tau_{h,j}}}}{\sqrt{ \abs{\mathcal{K}_i}}} \leq 1 $ indicates 
\[ \norm{G_h v_h}_{L^2(\tau_{h,j})} \lesssim \norm{ \nabla_{g_h} v_h}_{L^2(\mathcal{K}_i)} \]
which completes the proof.
\end{proof}

The boundedness of the operator $G^*_h$ in Algorithm \ref{alg:1} is a trivial case implicated from Lemma \ref{lem:boundedness}.  In the next we show the consistency of the PPPR gradient recovery operator by establishing the following lemma. 
\begin{lemma}
\label{lem:preserve}
Let $u\in  W^{3,\infty}(\manifold)$, and let $u_I$ be the linear interpolation of the function $u$ at every vertex of $\manifold_h$, then we have the estimate
\begin{equation}
\label{eq:preserve}
\norm{\nabla_g u - (T_h)^{-1} G_h u_I}_{0, \manifold} \le  h^2 \sqrt{\mathcal{A}(\manifold) }D(g,g^{-1}) \norm{u}_{3,\infty,\manifold},
\end{equation}
where $D(g,g^{-1})$ is a constant determined by the metric tensor $g$ and its inverse.
\end{lemma}
\begin{proof}
We start from a single triangle $\tau_{h,j}\subset \manifold_h$, and then go through all $j\in J_h$.
In particular, we consider the formulation \eqref{eq:local_gradient2} on each triangle.
Let $\tau_j \subset \manifold$ be the area corresponding to $\tau_{h,j}\subset \manifold_h$. Then we have 
\begin{eqnarray*}
\norm{\nabla_g u - (T_h)^{-1} G_h u_I}_{0,\tau_j}^2 
&= &\int_{\tau_{h,j}} \abs{\nabla \bar{u}_{\tau_{h,j}}(\partial \tilde{\f r}_{\tau_{h,j}})^\dag - G_h(u_I)}^2 \det(g\circ \tilde{\f r}_{\tau_{h,j}});\;
\end{eqnarray*}
where $ \bar{u}_{\tau_{h,j}}=u\circ \tilde{\f r}_{\tau_{h,j}}$ and $\tilde{\f r}_{\tau_{h,j}} $ is the geometric mapping from $\tau_{h,j}$ to $\tau_j$, that is: $\tau_j= \tilde{\f r}_{\tau_{h,j}}(\tau_{h,j}) $.
$G_h(u_I)|_{\tilde{\f r}_{\tau_{h,j}}}$ are the gradient values over the triangle $\tau_{h,j}$ by interpolating values recovered at the vertices of $\tau_{h,j}$. Therefore,  in the local coordinates of $\tau_{h,j}$, they are first order polynomials.

On the other hand, at every vertex $x_i$, let $\bar{u}_{i}(\zeta_i)=u\circ \f r_{h,i}(\zeta_i) = u\circ \f r_{i}(\zeta_i)$. The consistency of  polynomial preserving recovery operator $\bar{G}_h$ on the planar domain 
implies that 
\begin{equation} \label{equ:planar_con}
\abs{\nabla \bar{u}_{i}(\zeta_i) -\bar{G}_h \bar{u}_{I}(\zeta_i)}\leq Ch^2\norm{\bar{u}}_{3, \infty, \Omega_i} ,
\end{equation}
where $\zeta_i$ is the local coordinates for $x_i$. Let $\theta_i$ be coordinates for $x_i$ on $ \tau_{h,j}$.  Then we have 
\begin{equation}
\label{eq:local_equal}
\nabla_g u(x_i) =\nabla \bar{u}_{i}(\zeta_i) (\partial \f r_{i}(\zeta_i))^\dag = \nabla \bar{u}_{\tau_{h,j}}(\theta_i)(\partial \tilde{\f r}_{\tau_{h,j}}(\theta_i))^\dag, 
\end{equation}
and $G_h u_I(x_i)=\bar{G}_h \bar{u}_{I}(\zeta_i) (\bar{G}_h\f r_{h,i}(\zeta_i))^\dag $. 
Note that because both $\partial \f r_i$ and $ \bar{G}_h \partial \f r_{h,i}$ (see \eqref{eq:lowerbound}) are uniform bounded from below, using the consistency error estimation \eqref{eq:re_geo_error_bounds}, then we derive
\begin{equation}\label{equ:jjjj}
\abs{ (\partial \f r_i(\zeta_i))^\dag -  (\bar{G}_h \f r_{h,i} (\zeta_i))^\dag } \lesssim h^2\norm{\f r_i}_{3,\infty,\Omega_i}.
\end{equation}
By the triangle inequality and  the estimates  \eqref{equ:planar_con}  and \eqref{equ:jjjj}, we obtain 
\begin{equation*}
\begin{split}
 &\abs{\nabla \bar{u}_{i}(\zeta_i) (\partial \f r_{i}(\zeta_i))^\dag - G_h u_I(x_i)}\\
 \le & \abs{\nabla \bar{u}_{i}(\zeta_i) (\partial \f r_{i}(\zeta_i))^\dag - \bar{G}_h \bar{u}_{I}(\zeta_i) (\partial\f r_{i}(\zeta_i))^\dag} +
  \abs{\bar{G}_h \bar{u}_{I}(\zeta_i) (\partial\f r_{i}(\zeta_i))^\dag - \bar{G}_h \bar{u}_{I}(\zeta_i) (\bar{G}_h\f r_{h,i}(\zeta_i))^\dag}\\
    \lesssim & h^2\norm{\bar{u}}_{3, \infty, \Omega_i} \abs{(\partial\f r_{i}(\zeta_i))^\dag} + h^2\norm{\f r_i}_{3,\infty,\Omega_i}\abs{\bar{G}_h\bar{u}_I(\zeta_i)}\\
  \lesssim & h^2 \left(\norm{\bar{u}}_{3, \infty, \Omega_i}\norm{(\partial\f r_{i})^\dag}_{0,\infty,\Omega_i}  + \norm{\f r_i}_{3,\infty,\Omega_i}\norm{  \bar{u}_I}_{1,\infty,\Omega_i}\right).
\end{split}
\end{equation*}
Since $\nabla \bar{u}_{\tau_{h,j}}(\partial \tilde{\f r}_{\tau_{h,j}})^\dag$ is a vector valued function, then each of its components belongs to $W^{2,\infty}(\tau_{h,j})$.
Particularly, let $\left( \nabla \bar{u}_{\tau_{h,j}}(\partial \tilde{\f r}_{\tau_{h,j}})^\dag\right)_I$ be the linear interpolation of the vector-valued function $\nabla \bar{u}_{\tau_{h,j}}(\partial \tilde{\f r}_{\tau_{h,j}})^\dag$ on $\tau_{h,j}$. By the triangle inequality and the interpolating error estimate, we deduce that 
\begin{equation}
\label{eq:local_esti}
\begin{aligned}
& \norm{\nabla \bar{u}_{\tau_{h,j}}(\partial \tilde{\f r}_{\tau_{h,j}})^\dag - G_h(u_I)}_{0,\tau_{h,j}} \\
 \le & \norm{\nabla \bar{u}_{\tau_{h,j}}(\partial \tilde{\f r}_{\tau_{h,j}})^\dag - \left( \nabla \bar{u}_{\tau_{h,j}}(\partial \tilde{\f r}_{\tau_{h,j}})^\dag\right)_I }_{0,\tau_{h,j}}  + \\
 &  \norm{\left( \nabla \bar{u}_{\tau_{h,j}}(\partial \tilde{\f r}_{\tau_{h,j}})^\dag\right)_I  - G_h(u_I)}_{0,\tau_{h,j}}\\
\lesssim & h^2 \abs{\nabla \bar{u}_{\tau_{h,j}}(\partial \tilde{\f r}_{\tau_{h,j}})^\dag }_{2,\tau_{h,j}}  + \sum_{i\in V_j}\abs{\nabla \bar{u}_{\tau_{h,j}}(\theta_i)(\partial \tilde{\f r}_{\tau_{h,j}}(\theta_i))^\dag - G_h u_I(x_i)}\sqrt{\mathcal{A}(\tau_{h,j})}\\
\lesssim & h^2 \left(\abs{\nabla \bar{u}_{\tau_{h,j}}(\partial \tilde{\f r}_{\tau_{h,j}})^\dag }_{2,\tau_{h,j}}  + \sum_{i\in V_j} \norm{\bar{u}}_{3, \infty, \Omega_i}\norm{(\partial\f r_{i})^\dag}_{0,\infty,\Omega_i} \sqrt{\mathcal{A}(\tau_{h,j})}\right) \\
& + h^2 \sum_{i\in V_j} \norm{\f r_i}_{3,\infty,\Omega_i}\norm{\bar{u}_I}_{1,\infty,\Omega_i}\sqrt{\mathcal{A}(\tau_{h,j})}\\
\lesssim & h^2 \sum_{i\in V_j}  \left( \norm{\bar{u}}_{3, \infty, \Omega_i}\norm{(\partial\f r_{i})^\dag}_{2,\infty,\Omega_i} +\norm{\f r_i}_{3,\infty,\Omega_i}\norm{\bar{u}}_{1,\infty,\Omega_i}\right) \sqrt{\mathcal{A}(\tau_{h,j})};
  \end{aligned}
\end{equation}
where we have used   the relation \eqref{eq:local_equal}, and the following 
\[ \abs{\nabla \bar{u}_{\tau_{h,j}}(\partial \tilde{\f r}_{\tau_{h,j}})^\dag }_{2,\tau_{h,j}}= \abs{\nabla \bar{u}(\partial \f r_i)^\dag}_{2,\tau_{h,j}}  \leq  \norm{\bar{u}}_{3, \infty, \Omega_i}\norm{(\partial\f r_{i})^\dag}_{2,\infty,\Omega_i} \sqrt{\mathcal{A}(\tau_{h,j})},\]
in the last inequality. 
On every local domain $\Omega_i$, because of the regular property of $\manifold$, we have the following facts:
a),  $g\circ \f r_i=\partial \f r_i(\partial \f r_i)^T$, then $\norm{\partial \f r_i}_{k,\infty,\Omega_i}$ and $\norm{(\partial \f r_i)^\dag}_{k,\infty,\Omega_i}$ for $k\in \set{0,1,2}$ on all $\Omega_i$ can be estimated by $\sqrt{\norm{g}_{k,\infty}}$ and $\sqrt{\norm{g^{-1}}_{k,\infty}}$ respectively;
b),  Since $\manifold$ is $C^3$ smooth and regular with bounded curvature, therefore both $g$ and $g^{-1}$ and their derivatives up to second order are uniformly bounded from below and above.
On the other hand, we can estimate the norms
\[\norm{\bar{u}}_{k+1, \infty, \Omega_i}\leq \norm{(\det{g})^{-1}}_{0,\infty}\sqrt{\norm{g}_{k,\infty}}\norm{u}_{k+1,\infty,\manifold}  \text{ for } k\in \set{0,1,2}.\]
These allow us to have the estimate:
\begin{equation}
\label{eq:norm_equal}
\begin{aligned}
&\left(\norm{\bar{u}}_{3, \infty, \Omega_i}\norm{(\partial\f r_{i})^\dag}_{2,\infty,\Omega_i} +\norm{\f r_i}_{3,\infty,\Omega_i}\norm{\bar{u}}_{1,\infty,\Omega_i}\right)\\
\le & \sqrt{\norm{(\det{g})^{-1}}_{0,\infty}}\left(\norm{u}_{3, \infty, \manifold}\sqrt{\norm{g}_{2,\infty}}\sqrt{\norm{g^{-1}}_{2,\infty} }+\sqrt{\norm{g}_{2,\infty}}\sqrt{\norm{g}_{0,\infty}}\norm{u}_{1,\infty,\manifold}\right)\\
\le & C(g,g^{-1}) \norm{u}_{3, \infty, \manifold} ,
\end{aligned}
\end{equation}
where $C(g,g^{-1})$ are constants determined by the geometry of $\manifold$, and they are uniformly bounded whenever Assumption \ref{ass:irregular} satisfied.
Using the local estimate \eqref{eq:local_esti} and \eqref{eq:norm_equal}, we can further deduce that 
\begin{equation}
\label{eq:estimate1}
\begin{aligned}
& \int_{\tau_{h,j}} \abs{\nabla \bar{u}_{\tau_{h,j}}(\partial \tilde{\f r}_{\tau_{h,j}})^\dag - G_h(u_I)}^2 \det(g\circ \tilde{\f r}_{\tau_{h,j}})\;  \\
 \leq & \norm{\det g}_{0,\infty} \norm{\nabla \bar{u}_{\tau_{h,j}}(\partial \tilde{\f r}_{\tau_{h,j}})^\dag - G_h(u_I)}_{0,\tau_{h,j}}^2  \\
 \le & h^4 \abs{V_j}^2 \norm{\det g}_{0,\infty} (C(g,g^{-1}) \norm{u}_{3, \infty, \manifold} )^2\mathcal{A}(\tau_{h,j}).
\end{aligned}
\end{equation}
Note that $\abs{V_j}\equiv 3$ in our case.
Summing over both sides of \eqref{eq:estimate1} for all index $j\in J_h$, and taking the square root we get the final conclusion.
The constant
\[D(g,g^{-1})= \abs{V_j} \sqrt{\norm{\det g}_{0,\infty}} C(g,g^{-1})\]
where $C(g,g^{-1})$ is given in \eqref{eq:norm_equal}. 
Note that here the summation is bounded as we consider $\manifold$ to be compact and the fact that $\mathcal{A}(\manifold_h)\leq \mathcal{A}(\manifold)$, thus it does not reduce the order of $h$.
\end{proof}

Now we are ready to show the superconvergence of the recovered gradient  on $\manifold_h$. 

\begin{theorem}
\label{thm:superconvergence_error_surface}
Let Assumption \ref{ass:irregular} hold, and 
$u\in  W^{3,\infty}(\manifold)$ be the solution of \eqref{eq:var_laplace}, and $u_h$ be the solution of \eqref{eq:d_var_laplace}. Then
\begin{equation}
\label{eq:superconvergence_error_surface}
\begin{array}{ll}
 \norm{\nabla_g u - T_h^{-1}G_h u_h}_{0,\manifold} \le & h^2\left( \sqrt{\mathcal{A}(\manifold)}  D(g,g^{-1})\norm{u}_{3,\infty,\manifold}  + \norm{f}_{0,\manifold} \right)  + \\
 & C h^{1+\min \set{1,\sigma}}\left(\norm{u}_{3,\manifold}+\norm{u}_{2,\infty,\manifold} \right).
\end{array}
\end{equation}
where $D(g,g^{-1})$ is the same constant as Lemma \ref{lem:preserve}.
\end{theorem} 
\begin{proof}
This is readily shown by considering the triangle inequality
\begin{equation*}
\begin{array}{l}
 \norm{\nabla_g u - T_h^{-1}G_h u_h}_{0,\manifold}  \leq  \norm{\nabla_g u  - T_h^{-1}G_h u_I}_{0,\manifold} + \norm{T_h^{-1}G_h ( u_I- u_h)}_{0,\manifold}.
\end{array}
\end{equation*}
The first term is bounded by  Lemma \ref{lem:preserve}. 
For the second term, since both $(T_h)^{-1} $ and $G_h $ are bounded operators ( Lemma \ref{lem:transform} and Lemma \ref{lem:preserve}), we have
\[\norm{(T_h)^{-1} G_h  (u_I -  u_h)}_{0,\manifold} \leq  C\norm{\nabla_{g_h} (u_I - u_h) }_{0,\manifold_h}\]
then using the result
\footnote{The $\mathcal{O}(h^{2\sigma})$ condition is asked for the projected triangle meshes on each $\Omega_i$ in order to show the supercloseness, while what we have assumed is in fact on the meshes before projection as \cite{WeiChenHuang2010}.
We argue that for general smooth surfaces with uniformly bounded curvature, using the ways described in our algorithms, the projected shape of meshes will not be significantly changed as \cite{WeiChenHuang2010}, therefore the  $\mathcal{O}(h^{2\sigma})$ condition can be guaranteed, although this might be not the case for the meshes located at the high curvature areas.
Once a surface is highly curved, one may have to take into account the ratio of the high curvature areas, thus $\mathcal{O}(h^{2\sigma})$ condition may be adapted to a new index $\bar{\sigma}$ according to ratio of the high curvature areas.
But in this paper, we skip the quantitative discussion on this point. } 
of \cite[Theorem 3.5]{WeiChenHuang2010} to estimate $\norm{\nabla_{g_h} (u_I- u_h) }_{0,\manifold_h}$.
These lead to the final estimate.
\end{proof}

Due to Lemma \ref{lem:transform}, we have the following result immediately, which is verified in our numerical part Section \ref{sec:numerics}. 
\begin{corollary}
\label{cor:superconvergence}
Let the same assumptions as Theorem \ref{thm:superconvergence_error_surface} hold.  Then
\begin{equation}
\label{eq:superconvergence2}
\begin{array}{ll}
\norm{T_h\nabla_g u - G_h u_h}_{0,\manifold_h} \leq &  h^2\left( \sqrt{\mathcal{A}(\manifold)}  D(g,g^{-1})\norm{u}_{3,\infty,\manifold}  + \norm{f}_{0,\manifold} \right) + \\
&  Ch^{1+\min \set{1,\sigma}}\left(\norm{u}_{3,\manifold}+\norm{u}_{2,\infty,\manifold} \right).
\end{array}
\end{equation}
\end{corollary} 
\section{Recovery-based a posteriori error estimator}
\label{sec:estimator}

The gradient recovery operator $G_h$ naturally provides an {\it a posteriori } error estimator.
We define a local {\it a posteriori} error estimator on each triangular element $ \tau_{h,j} $ as
\begin{equation}\label{equ:localind}
\eta_{h,\mathcal{\tau}_{h,j}} =
\|G_hu_h - \nabla_{g_h} u_h\|_{0, \tau_{h,j}}, 
\end{equation}
and the corresponding global error estimator  as
\begin{equation}\label{equ:globalind}
\eta_h = \left( \sum\limits_{j\in J_h}\eta_{h,\tau_{h,j}}^2\right)^{1/2}.
\end{equation}

With the previous superconvergence result, we can  show the asymptotic exactness of error estimators based on the recovery operator $G_h$.

\begin{corollary}\label{cor:asyexact}
Assume the same conditions in Theorem \ref{thm:superconvergence_error_surface} and let $u_h$ be the  finite element solution of discrete variational problem  \eqref{eq:d_var_laplace}.
Further assume that there is a constant $C(u)>0$ such that
 \begin{equation}\label{equ:satassum}
\norm{T_h\nabla_g u-\nabla_{g_h} u_h}_{0,\manifold_h} \ge C(u) h.
\end{equation}
Then it holds that
\begin{equation}
\left| \frac{\eta_h}{\norm{ T_h \nabla_{g} u-\nabla_{g_h} u_h}_{0,\manifold_h}}  -1 \right| \lesssim h^{\min\set{1,\sigma}}.
\end{equation}
\end{corollary}
\begin{proof}
By the triangle inequality,  we have 
\begin{equation*}
\eta_h \le \|G_hu_h - T_h \nabla_{g} u\|_{{0,\manifold_h}} +  \| T_h \nabla_{g} u-\nabla_{g_h} u_h\|_{0,\manifold_h}
\end{equation*}
and hence 
\begin{equation*}
\left|\frac{\eta_h}{\norm{ T_h \nabla_{g} u-\nabla_{g_h} u_h}_{0,\manifold_h}}  -1 \right| \le \frac{\|G_hu_h - \nabla_{g_h} u_h\|_{{0,\manifold_h}}}{ \| T_h \nabla_{g} u-\nabla_{g_h} u_h\|_{0,\manifold_h}} \lesssim h^{\min\set{1,\sigma}}.
\end{equation*}
where we use the superconvergence result \eqref{eq:superconvergence2} and the  assumption \eqref{equ:satassum} 
in the last inequality.
\end{proof}

\begin{remark}
The assumption \eqref{equ:satassum} is common   assumption to show the asymptotical exactness of recovery-based  {\it a posteriori} error estimators as \cite{AinsworthOden2000, NagaZhang2004, ZhangNaga2005}.  
It is reasonable since  that the finite element solution error is not better than the interpolation error  which is bounded from below by $\mathcal{O}(h)$(except some trivial cases).
\end{remark}

\begin{remark}
Corollary \ref{cor:asyexact} implies that \eqref{equ:localind} (or \eqref{equ:globalind}) is an asymptotically exact {\it a posteriori } error estimator for surface finite element methods.
\end{remark}

\section{Numerical Results}
\label{sec:numerics}

In this section, we present several numerical examples to demonstrate the superconvergence property of the proposed gradient recovery operators and make comparisons with existing gradient recovery operators.  The first example is to show the superconvergence results of the proposed gradient recovery operators even though the element patch is not $\mathcal{O}(h^2)$-symmetric. In this example, the vertices are located exactly on the torus.
The second one is to compare the results on a more complicated surface and to demonstrate the superiority of the PPPR method for surfaces with high curvature, which is an example that the vertices are not located on the exact surfaces.
The last two are to show the asymptotic exactness of the recovery-based {\it a posterior } error estimator introduced in Section \ref{sec:estimator}. Some of our numerical tests are conducted based on the  MATLAB package {\emph i}FEM \cite{chen2009ifem}.  Except for the first example, the initial meshes for the other three examples are generated using the three-dimensional surface mesh generation module of the Computational Geometry Algorithms Library \cite{cgal}. 
To get meshes in other levels, we first perform either the uniform refinement or the newest bisection \cite{Chen2008}. Then we project the newest vertices onto the $\manifold$. 
In the general case, there is no explicit project map available.  Hence we adopt the first order approximation of projection map as given in \cite{DemlowDziuk2007}. 
Thus, the vertices of the meshes are not on the exact surface $\manifold$ but in an $h^2$ neighbourhood for the second and fourth example.
We notice that in such cases the superconvergence results can still be observed.

Let $G_h^{SA}$, $G_h^{WA}$, and $G_h^{ZZ}$ be recovery operators by simple averaging, weighted averaging, and Zienkiewicz-Zhu schemes on tangent planes \cite{WeiChenHuang2010}, respectively. 
Note  that we use the exact normal vectors for $G_h^{ZZ}$ in the numerical examples.
We denote $G^*_{h}$, $G_{h}$, and $G_h^a$ to be the recovery operators given by Algorithm \ref{alg:1}, Algorithm \ref{alg:2} and Algorithm \ref{alg:1} with approximations of normal vectors, respectively.
The approximating normal vectors are computed by weighted averaging for the tests with $G^a_h$ in our examples, which are also used to implement Algorithm \ref{alg:2} to construct the local parametric domains $\Omega_i$.
Another remark is that we use the function value preserving skill for the PPPR $G_h$, but not for $G^*_h$.
For the reason of making comparisons, we define:
\begin{align*}
 &  De=\norm{T_h\nabla_g u -\nabla_{g_h} u_h }_{0, \manifold_h}, &
 &  De^I = \norm{\nabla_{g_h} u_I - \nabla_{g_h} u_h}_{0, \manifold_h},\\
 &  De^{r_1} = \norm{T_h\nabla_g u - G^*_{h} u_h}_{0, \manifold_h},&
 &  De^{r_2} = \norm{T_h\nabla_g u - G_{h} u_h}_{0,  \manifold_h},\\
 &  De^{r_3} = \norm{T_h\nabla_g u - G_h^a u_h}_{0,\manifold_h},&
 &  De^{SA} = \norm{T_h\nabla_g u - G_h^{SA} u_h}_{0, \manifold_h}, \\
 &  De^{WA} = \norm{T_h\nabla_g u - G_h^{WA} u_h}_{0, \manifold_h},&
 &  De^{ZZ} = \norm{T_h\nabla_g u - G_h^{ZZ} u_h}_{0, \manifold_h};
\end{align*}
where $u_h$ is the finite element solution, $u$ is the analytical solution and $u_I$ is the linear finite element interpolation of $u$.

In Numerical Example $2$, we shall compare the discrete maximal errors of the above six discrete gradient recovery methods. For that reason, we introduce the following notations
\begin{align*}
 &  De_0^{r_1} = \norm{T_h\nabla_g u - G^*_{h} u_h}_{0, \infty, \manifold_h},&
 &  De_0^{r_2} = \norm{T_h\nabla_g u - G_{h} u_h}_{0, \infty, \manifold_h},\\
 &  De_0^{r_3} = \norm{T_h\nabla_g u - G_h^a u_h}_{0, \infty, \manifold_h},&
 &  De_0^{SA} = \norm{T_h\nabla_g u - G_h^{SA} u_h}_{0, \infty, \manifold_h}, \\
 &  De_0^{WA} = \norm{T_h\nabla_g u - G_h^{WA} u_h}_{0, \infty,\manifold_h},&
 &  De_0^{ZZ} = \norm{T_h\nabla_g u - G_h^{ZZ} u_h}_{0, \infty, \manifold_h};
\end{align*}
where $\|\cdot\|_{0, \infty, \manifold_h}$ means the maximum absolute value at all vertices.

In the following tables, all convergence rates are listed in term of the degree of freedom(Dof).
Noticing  $\text{Dof}  \approx h^{-2}$, the corresponding convergence rates in term of the mesh size $h$ are double of what we present in the tables.

\subsection{Numerical Example 1} Our first example is to consider Laplace-Beltrami equation on a torus surface. The right hand function $f$ is chosen to fit the exact solution $u(x,y,z) = x-y$.
The  signed 
distance function of torus surface is 
\begin{equation}
\Phi(x) = \sqrt{(4-\sqrt{x_1^2+x_2^2)^2+x_3^2}}-1.
\end{equation}

\begin{figure}[h!]
  \begin{minipage}[b]{.45\linewidth}
   \centering
   \includegraphics[width=0.8\textwidth]{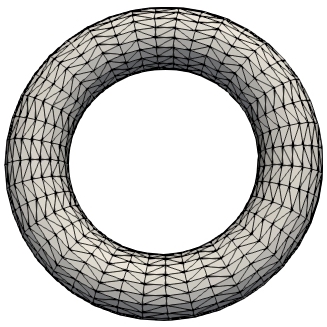}
  \subcaption{}\label{fig:1a}
   \end{minipage}%
 \begin{minipage}[b]{.45\linewidth}
  \centering
 \includegraphics[width=0.8\textwidth]{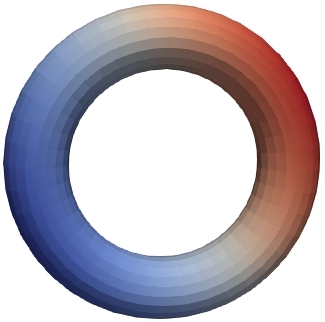}
  \subcaption{}\label{fig:1b}
 \end{minipage}
 \caption{Numerical Solution on Torus Surface: (a) Mesh; (b) Solution.} 
 \label{fig:torus}
\end{figure}

To construct a series meshes on torus without $\mathcal{O}(h^2)$ symmetric property of their element patches, 
we firstly make a series of uniform meshes of Chevron pattern and map the mesh onto the torus.
Figure  \ref{fig:torus} plots  the uniform mesh with 800 Dofs  and  the corresponding finite element solution.

   Table  \ref{tab:torus} lists the numerical results.  As expected,  $H^1$ error 
   of the  finite element solution is  of $\mathcal{O}(h)$.  Since the generated uniform meshes
  satisfy the $\mathcal{O}(h^{2\sigma})$ condition,    $\mathcal{O}(h^2)$ supercloseness
  for $De^I$ is observed.  Concerning the convergence of recovered gradients, both the recovered gradient by PPR with exact normal field and by the PPPR  have a superconvergence rate of order $\mathcal{O}(h^2)$;  while the recovered gradient using PPR with approximated normal field and the other three methods in \cite{WeiChenHuang2010} only converge at the optimal rate $\mathcal{O}(h)$.

\vspace{-1cm}
\begin{table}[h!]
\centering
\caption{Numerical Results for equation \eqref{eq:var_laplace} on torus surface.}\label{tab:torus}
\small 
\begin{tabular}{|c|c|c|c|c|c|c|c|c|c|c|c|}
\hline 
Dof& $De$&order&$De^I$&order&$De^{r_1}$&Order&$De^{r_2}$&order\\ \hline
200&2.52e+00&--&9.43e-01&--&1.50e+00&--&1.59e+00&--\\ \hline
800&1.26e+00&0.50&2.65e-01&0.92&4.12e-01&0.93&4.37e-01&0.93\\ \hline
3200&6.29e-01&0.50&6.92e-02&0.97&1.06e-01&0.98&1.13e-01&0.98\\ \hline
12800&3.14e-01&0.50&1.75e-02&0.99&2.67e-02&0.99&2.84e-02&0.99\\ \hline
51200&1.57e-01&0.50&4.40e-03&1.00&6.70e-03&1.00&7.12e-03&1.00\\ \hline
204800&7.86e-02&0.50&1.10e-03&1.00&1.67e-03&1.00&1.78e-03&1.00\\ \hline
819200&3.93e-02&0.50&2.75e-04&1.00&4.19e-04&1.00&4.45e-04&1.00\\ \hline
3276800&1.97e-02&0.50&6.88e-05&1.00&1.05e-04&1.00&1.11e-04&1.00\\ \hline
\hline 
Dof&$De^{r_3}$&order&$De^{SA}$&order&$De^{WA}$&Order&$De^{ZZ}$&order\\ \hline
200&1.52e+00&--&2.27e+00&--&2.28e+00&--&2.27e+00&--\\ \hline
800&4.74e-01&0.84&7.22e-01&0.83&7.25e-01&0.83&6.91e-01&0.86\\ \hline
3200&1.68e-01&0.75&2.48e-01&0.77&2.49e-01&0.77&2.19e-01&0.83\\ \hline
12800&7.18e-02&0.61&1.03e-01&0.63&1.03e-01&0.63&8.39e-02&0.69\\ \hline
51200&3.42e-02&0.54&4.86e-02&0.54&4.86e-02&0.54&3.80e-02&0.57\\ \hline
204800&1.69e-02&0.51&2.39e-02&0.51&2.39e-02&0.51&1.84e-02&0.52\\ \hline
819200&8.40e-03&0.50&1.19e-02&0.50&1.19e-02&0.50&9.16e-03&0.51\\ \hline
3276800&4.20e-03&0.50&5.94e-03&0.50&5.94e-03&0.50&4.57e-03&0.50\\ \hline
\end{tabular}
\end{table}

\subsection{Numerical Example 2}
In this example, we take a surface \cite{DziukElliott2013} which contains high curvature features.
It can be represented as the zero level of the following level set function 
\begin{equation*}
\Phi(x) = \frac{1}{4}x_1^2+x_2^2+\frac{4x_3^2}{(1+\frac{1}{2}\sin(\pi x_1))^2}-1.
\end{equation*}
We consider the Laplace-Beltrami equation \eqref{eq:laplace} with exact solution $u = x_1x_2$.
The right-hand side function $f$ can be computed from $u$. 

\begin{figure}
   \centering
   \subcaptionbox{Mesh\label{fig:elliotmesh}}
  {\includegraphics[width=0.68\textwidth]{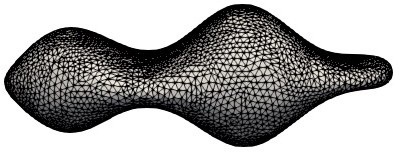}}
  \subcaptionbox{Solution\label{fig:elliotsol}}
   {\includegraphics[width=0.75\textwidth]{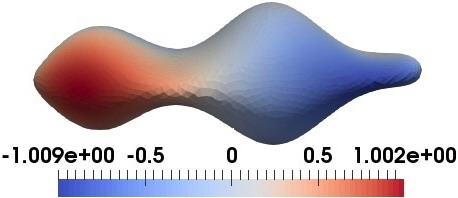}}
  \subcaptionbox{Error distribution of recovered gradient\label{fig:elliotdist}}
  {\includegraphics[width=0.75\textwidth]{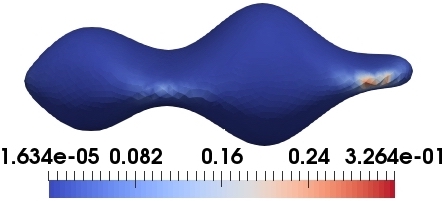}}
   \caption{Numerical Solution on general surface.}\label{fig:elliot}
\end{figure}

Figure \ref{fig:elliotsol} shows the finite element solution $u_h$ on Delaunay mesh, see \ref{fig:elliotmesh},  with 4606 Dofs. 
The numerical results is reported in Table \ref{tab:elliot}. From the table, we clearly see that 
$De$ converges at the optimal rate $\mathcal{O}(h)$ and $De^I$ converges at a superconvergent rate $\mathcal{O}(h^2)$.
As demonstrated in \cite{Chernyshenk02015}, some regions of the surface show significant high curvature. 
Due to the existence of these areas, only sub-superconvergence rate of order $\mathcal{O}(h^{1.8})$ is observed for 
PPR with approximated normal field and the other three methods in \cite{WeiChenHuang2010}. 
In contrast, the $\mathcal{O}(h^2)$ superconvergence rate can be observed in the PPR with exact normal field and in the PPPR method.
To look more clearly into the relations between the recovery accuracy and the high curvature of a surface, we add another set of comparison in this example.
In our numerical tests, we observed that the maximal recovery errors always happened in the area of the meshes generated from highest curvature surface regions.
We plot a case example of the distribution of the error function $|G_h u_h- T_h\nabla_g u|$ in Figure \ref{fig:elliotdist}.
Table \ref{tab:elliotl0} reports the maximal discrete errors of all the above six gradient recovery methods,
in which  PPPR method is the only one to achieve the superconvergence rate of $\mathcal{O}(h^2)$ asymptotically in the discrete maximal norm.
This gives the evidence to our statement in Remark \ref{rem:highcurvature} that PPPR is relatively \emph{curvature stable} compared to the other methods. 
At that point, we can say that PPPR is the best one for arbitrary meshes and meshes generated by high curvature surfaces.
Thus, in the following two examples, we shall only consider the PPPR method. 

In this  example, we find from Table \ref{tab:elliot} and \ref{tab:elliotl0} that the results of the Algorithm \ref{alg:1} ($De^{r_1}$),
which uses the exact normal vectors, are worse than the results of Algorithm \ref{alg:2} ($De^{r_2}$).
This is not surprising, as we have reported that in this complicated surface case, 
the vertices of the discrete mesh are not located on the exact analytical surface any more. 
Therefore even with exact normal vectors, it brings unavoidable errors to the computations.
This also shows an advantage of the PPPR method (Algorithm \ref{alg:2}).

\begin{table}[htb!]
\centering
\caption{Numerical Results for equation \eqref{eq:var_laplace} on a general surface}\label{tab:elliot}
\small
\begin{tabular}{|c|c|c|c|c|c|c|c|c|c|c|}
\hline 
Dof& $De$&order&$De^I$&order&$De^{r_1}$&Order&$De^{r_2}$&order\\ \hline
1153&5.46e-01&--&2.78e-01&--&4.77e-01&--&3.34e-01&--\\ \hline
4606&2.85e-01&0.47&1.18e-01&0.62&2.01e-01&0.62&1.29e-01&0.69\\ \hline
18418&1.40e-01&0.51&3.45e-02&0.89&6.58e-02&0.81&4.38e-02&0.78\\ \hline
73666&6.97e-02&0.50&9.86e-03&0.90&1.97e-02&0.87&1.28e-02&0.89\\ \hline
294658&3.48e-02&0.50&2.58e-03&0.97&5.33e-03&0.95&3.40e-03&0.96\\ \hline
1178626&1.74e-02&0.50&6.57e-04&0.99&1.37e-03&0.98&8.68e-04&0.99\\ \hline
4714498&8.70e-03&0.50&1.66e-04&0.99&3.46e-04&0.99&2.18e-04&1.00\\ \hline
\hline
Dof&$De^{r_3}$&order&$De^{SA}$&order&$De^{WA}$&Order&$De^{ZZ}$&order\\ \hline
1153&4.71e-01&--&4.83e-01&--&4.86e-01&--&4.95e-01&--\\ \hline
4606&1.98e-01&0.62&2.26e-01&0.55&2.30e-01&0.54&2.18e-01&0.59\\ \hline
18418&6.63e-02&0.79&8.30e-02&0.72&8.59e-02&0.71&7.45e-02&0.78\\ \hline
73666&2.06e-02&0.84&2.69e-02&0.81&2.82e-02&0.80&2.33e-02&0.84\\ \hline
294658&5.87e-03&0.91&7.72e-03&0.90&8.27e-03&0.89&6.60e-03&0.91\\ \hline
1178626&1.64e-03&0.92&2.14e-03&0.93&2.36e-03&0.90&1.83e-03&0.93\\ \hline
4714498&4.70e-04&0.90&6.04e-04&0.91&6.97e-04&0.88&5.22e-04&0.90\\ \hline
\end{tabular}
\end{table}

\begin{table}[htb!]
\centering
\caption{Comparison of discrete maximal norms of gradient recovery methods on a general surface}\label{tab:elliotl0}
\begin{tabular}{|c|c|c|c|c|c|c|c|c|c|c|}
\hline 
Dof& $De_0^{r_1}$&order&$De_0^{r_2}$&order&$De_0^{r_3}$&Order\\ \hline
1153&9.70e-01&--&7.93e-01&--&8.73e-01&--\\ \hline
4606&5.43e-01&0.42&3.26e-01&0.64&4.77e-01&0.44\\ \hline
18418&1.92e-01&0.75&1.09e-01&0.79&2.22e-01&0.55\\ \hline
73666&8.57e-02&0.58&5.18e-02&0.54&9.16e-02&0.64\\ \hline
294658&2.50e-02&0.89&1.40e-02&0.94&3.51e-02&0.69\\ \hline
1178626&7.80e-03&0.84&3.59e-03&0.98&1.59e-02&0.57\\ \hline
4714498&3.56e-03&0.57&9.03e-04&0.99&7.57e-03&0.53\\ \hline
\hline
Dof&$De_0^{SA}$&order&$De_0^{WA}$&order&$De_0^{ZZ}$&Order\\ \hline
1153&7.16e-01&--&7.09e-01&--&8.13e-01&--\\ \hline
4606&5.09e-01&0.25&5.36e-01&0.20&5.83e-01&0.24\\ \hline
18418&2.73e-01&0.45&3.05e-01&0.41&2.65e-01&0.57\\ \hline
73666&1.42e-01&0.47&1.47e-01&0.53&1.11e-01&0.63\\ \hline
294658&5.66e-02&0.66&6.08e-02&0.64&3.67e-02&0.80\\ \hline
1178626&2.39e-02&0.62&2.75e-02&0.57&1.62e-02&0.59\\ \hline
4714498&1.12e-02&0.55&1.31e-02&0.54&7.63e-03&0.54\\ \hline
\end{tabular}
\end{table}

\subsection{Numerical Example 3}
In the example, we consider a benchmark problem for adaptive finite element method for the Laplace-Beltrami equation on the sphere \cite{DemlowDziuk2007,Demlow2012,Chernyshenk02015}.  We choose the right-hand side function $f$ such that the exact solution in spherical coordinate is 
given by 
\begin{equation*}
u = \sin^{\lambda}(\theta)\sin(\psi).
\end{equation*}
In case of $\lambda <1$, it easy to see that the solution $u$ has two singularity points at north and south poles and the solution $u$ is barely in $H^1(\manifold)$.
In fact, $u \in H^{1+\lambda}(\manifold)$.

\begin{figure}
  \begin{minipage}[b]{.5\linewidth}
   \centering
 \includegraphics[width=0.8\textwidth]{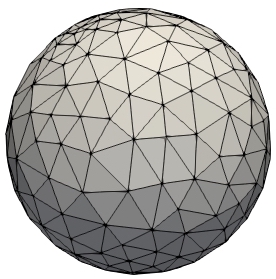}
  \subcaption{} \label{fig:sphere_init}
  \end{minipage}%
  \begin{minipage}[b]{.5\linewidth}
  \centering
\includegraphics[width=0.8\textwidth]{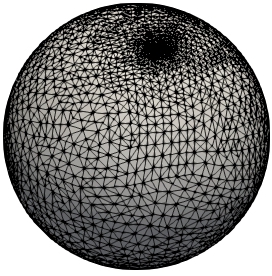}
   \subcaption{}\label{fig:sphere_adaptive}
  \end{minipage}
\caption{Meshes for Example 3:  (a)  Initial mesh; (b) Adaptively refined mesh.}
 \label{fig:sphere_mesh}
  \begin{minipage}[b]{.5\linewidth}
    \centering
  \includegraphics[width=0.8\textwidth]{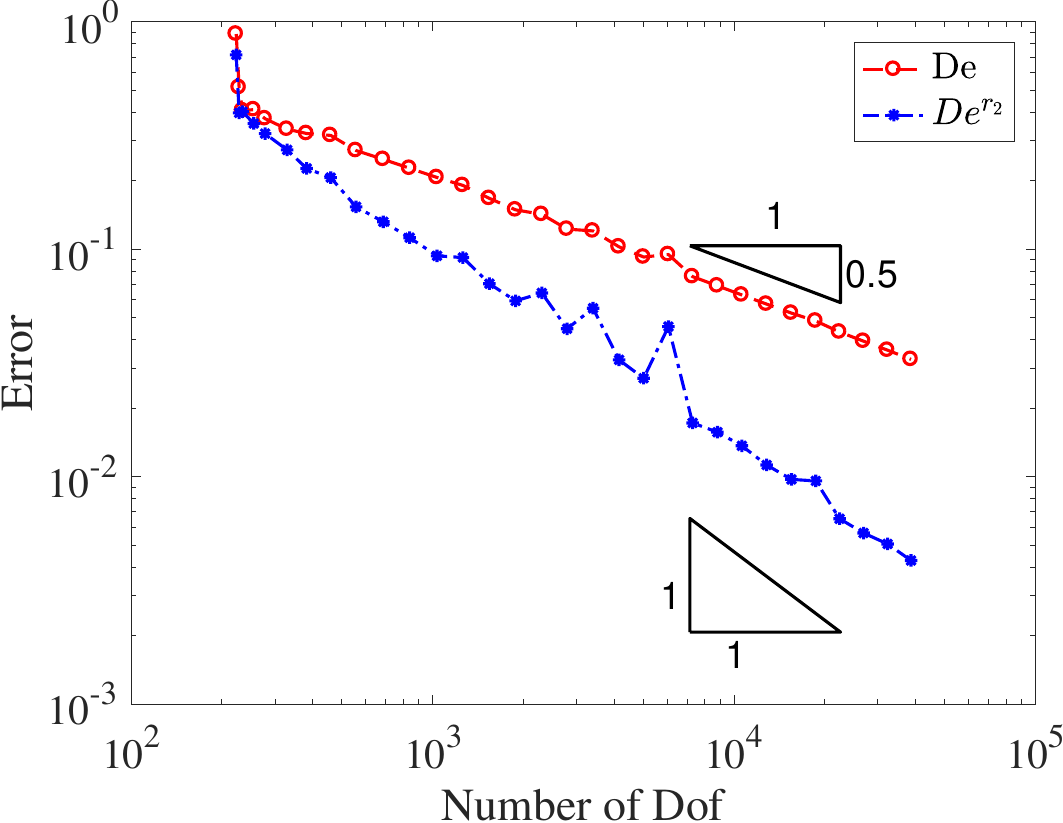}
   \subcaption{} \label{fig:sphere_err}
  \end{minipage}%
  \begin{minipage}[b]{.5\linewidth}
  \centering
\includegraphics[width=0.8\textwidth]{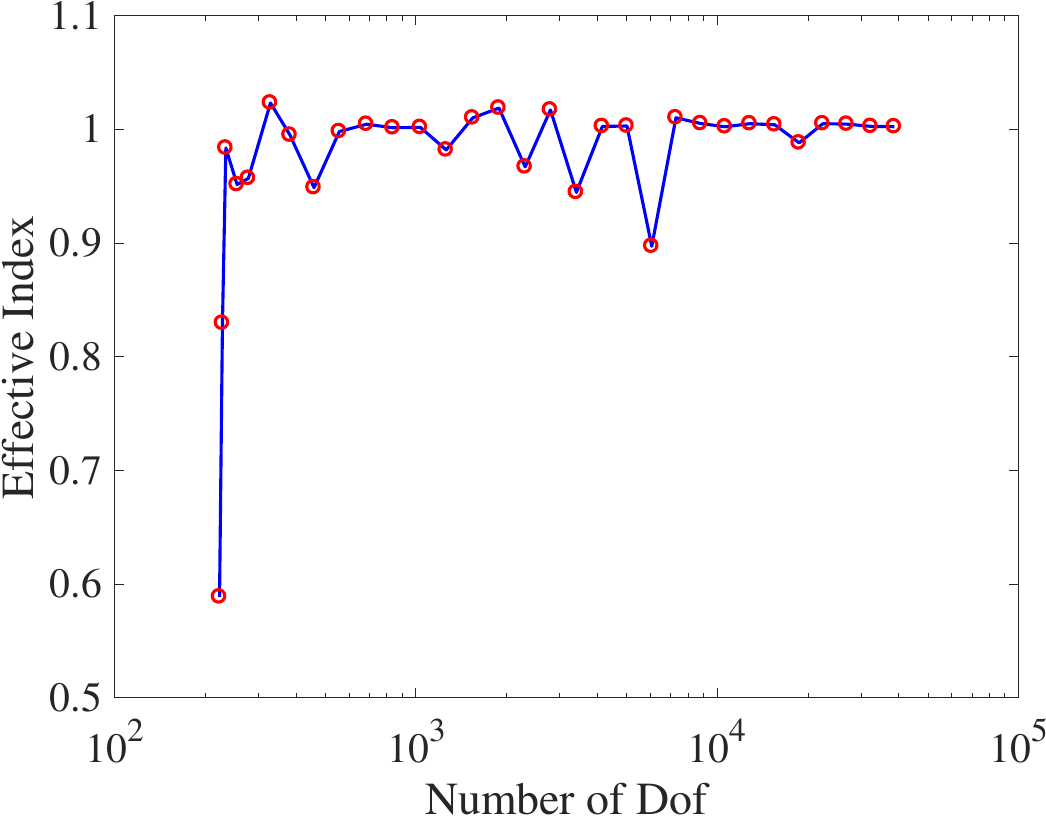}
  \subcaption{}\label{fig:sphere_idx}
  \end{minipage}
\caption{Numerical Result for Example 3:  (a)  Errors; (b) Effective index.}
  \label{fig:sphere_result}
\end{figure}

To obtain the optimal convergence rate, we use the adaptive finite element method (AFEM).  Different from the existing methods in the literature, the recovery-based \textit{a posteriori} error estimator is adopted.  We start with the initial mesh given as in Fig \ref{fig:sphere_init}. 
The mesh is adaptively refined using the Do\"rfler \cite{Dorfler1996} marking strategy with parameter equal to  $0.3$.
Fig \ref{fig:sphere_adaptive} plots the mesh after the 18 adaptive refinement steps. The mesh successfully resolves the singularities. 
The numerical errors are displayed in Fig \ref{fig:sphere_err}. As expected, an optimal convergence rate for $H^1$ error can be observed. Also,  we observe that the recovered gradient is superconvergent to the exact gradient at a rate of $\mathcal{O}(h^2)$. 

To test  the performance of our new recovery-based \textit{a posterior} error estimator for the Laplace-Beltrami problem, 
the effectivity index $\kappa$ is used to measure the quality of
an error estimator \cite{AinsworthOden2000,Babuska2001}, which  is defined by the ratio between the estimated error and the  exact error
\begin{equation}
    \kappa = \frac{\norm{
    G_h u_h - \nabla_{g_h} u_h}_{0, \manifold_h}}
    {\norm{ T_h\nabla_g u - \nabla_{g_h}u_h)}_{0, \manifold_h}}
    \label{equ:effect}
\end{equation}
The effectivity index is plotted in  Fig \ref{fig:sphere_idx} .
We see that $\kappa$ converges asymptotically to $1$ 
which  indicates the posteriori error estimator 
\eqref{equ:localind} (or \eqref{equ:globalind} )
is asymptotically exact.

\subsection{Numerical Example 4}  
\begin{figure}[h!]
   \begin{minipage}[b]{.48\linewidth}
    \centering
   \includegraphics[width=0.8\textwidth]{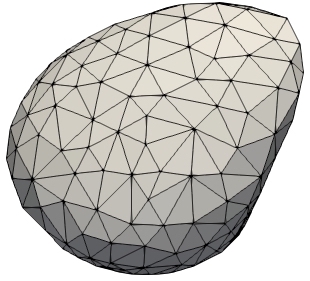}
   \subcaption{} \label{fig:dziuk_init}
   \end{minipage}%
   \begin{minipage}[b]{.48\linewidth}
    \centering
   \includegraphics[width=0.8\textwidth]{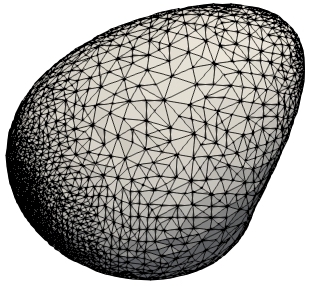}
   \subcaption{}\label{fig:dziuk_adaptive}
   \end{minipage}
\caption{Meshes for Example 4:  (a)  Initial mesh; (b) Adaptively refined mesh.}
  \label{fig:dziuk_mesh}
\end{figure}

\begin{figure}[h!]
     \begin{minipage}[b]{.48\linewidth}
    \centering
    \includegraphics[width=0.8\textwidth]{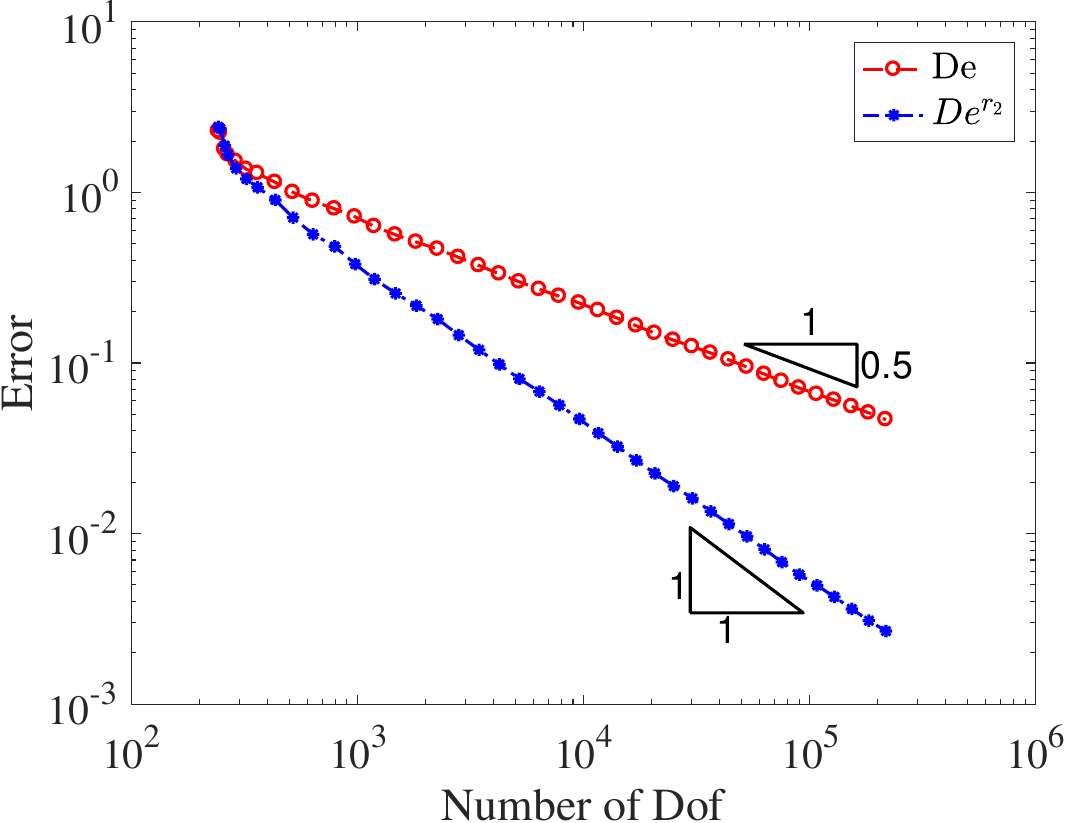}
    \subcaption{} \label{fig:dziuk_err}
    \end{minipage}%
    \begin{minipage}[b]{.48\linewidth}
    \centering
    \includegraphics[width=0.8\textwidth]{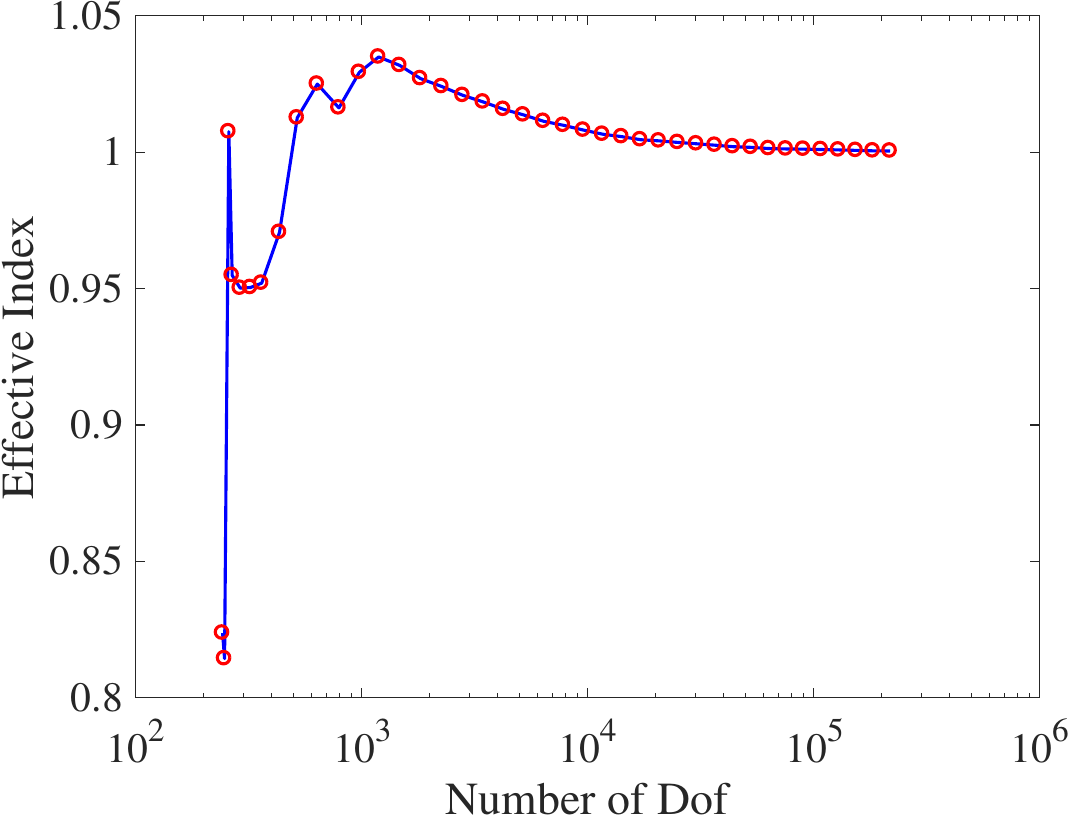}
    \subcaption{}\label{fig:dziuk_idx}
    \end{minipage}
\caption{Numerical Result for Example 4:  (a)  Errors; (b) Effective index.}
 \label{fig:dziuk_result}
\end{figure}

In this example, we consider the following  Laplace-Beltrami type equation on Dziuk surface as in \cite{DednerMadhavanStinner2016}:
\begin{equation*}
-\Delta_{g}+u = f, \quad \text{on } \Gamma,
\end{equation*}
where $\Gamma = \left\{x\in \mathbb{R}^3: (x_1-x_3^2)^2+x_2^2+x_3^2=1\right\}$.  $f$ is chosen to fit the exact solution
\begin{equation*}
u(x,y,z) = e^{\frac{1}{1.85-(x-0.2)^2}}\sin(y). 
\end{equation*}
Note that the solution has an exponential peak.  To track this phenomenon, we adopt AFEM with an initial mesh graphed in Fig \ref{fig:dziuk_init}. Fig \ref{fig:dziuk_adaptive} shows the adaptively refined mesh. We would like to point out that the mesh is refined not only around the exponential peak but also at the high curvature areas.
Fig \ref{fig:dziuk_err} displays the numerical errors. It demonstrates the optimal convergence rate in $H^1$ norm and a superconvergence rate for the recovered gradient.   The effective index is shown in Fig \ref{fig:dziuk_idx}, which converges to 1 quickly after the first few iterations. Again, it  indicates the error estimator  \eqref{equ:localind} (or \eqref{equ:globalind} ) is asymptotically exact.

\section{Conclusion}
In this paper, we have proposed a novel gradient recovery method which uses parametric polynomials to fit both manifolds and also FEM solutions defined on the manifolds, and then to recover a better gradient for the FEM solutions. 
In comparing with existing methods for data on surfaces in the literature, cf. \cite{DuJu2005,WeiChenHuang2010},
the proposed method has several improvements:
The first highlight is that it does not requires the normal vectors of the exact manifold, which makes it a realistic and robust method for practical problems;
Second, it does not need the element patch to be $\mathcal{O}(h^2)$ symmetric to achieve superconvergence;
Third, all of our numerical tests show evidence that it is a curvature stable method in comparing with the existing methods.
We have advanced the traditional PPR method (for planar problems) to function value preserving in the meantime,
and shown the capability of the recovery operator for constructing {\it a posteriori} error estimator.
In fact, the superconvergence does hold for $\manifold_h$ with no exact vertices, but the theoretical proof is postponed to a follow-up paper. 
Even though we have only discussed linear finite element methods on triangulated meshes,
the idea should be applicable to higher order FEM on more accurate approximations of surfaces, e.g., piecewise polynomial surfaces, B-splines or NURBS.
However, these are non-trivial works, and we leave them as future topics.

Gradient recovery has other applications, like enhancing eigenvalues\cite{NagaZhangZhou2006, NagaZhang2012, GuoZhangZhao2017}, simplifying higher order discretization of PDEs\cite{GrandeReusken2016}, designing new numerical methods for higher order PDEs\cite{ChenGuoZhangZou2017,GuoZhangZou2017, GZZ2018b, XGZ2019}. 
Moreover, it may help for the vector field regularization in the context of \cite{DonJueSchTak17}, where the geometric approximation accuracy is asked to be one more extra order higher than the order of function approximation accuracy in order to optimally regularizing vector fields on manifolds.
The superconvergence property of the recovery scheme might be able to reduce the additional higher order requirement in manifolds approximation for vector fields, and achieve optimal convergence rates as the case of scalar valued functions.
It would be interesting to investigate further the full usage of the PPPR method for problems with solutions defined on manifolds.

\section*{Ackowledgement}
The authors thank Dr. Pravin Madhavan,  Dr. Bjorn Stinner and Dr. Andreas Dedner for their kind help and discussions on a numerical example. 
The authors also thank the editor and referees for their efforts and valuable comments and suggestions which have significantly improved the paper.

\appendix
\section{Proof of Lemma \ref{lem:invariant}}
\label{appendix}
\begin{proof}
In general, there are infinitely many isomorphic parameterizations for a given patch $S\subset \manifold$.
Let us pick arbitrarily two of them, which are denoted by
\begin{eqnarray*}
   \f r: \Omega \rightarrow  S \quad  \text{ and } \quad \f s: \Omega_s \rightarrow  S \;, 
    \end{eqnarray*} 
respectively, where $\Omega$ and $\Omega_s$ are planar parameter domains, then there exist
   \begin{eqnarray*}
      \f t : \Omega &\rightarrow & \Omega_s
   \end{eqnarray*} 
to be a bijective, differentiable mapping, such that $\f r =\f s\circ \f t$.
That means for an arbitrary but fixed position $x\in S$, we have $\xi\in \Omega $ and $\f t(\xi)=\zeta$, such that 
  \begin{equation*}
x=  \f s (\zeta)  = \f s( \f t(\xi))=\f r(\xi) .  
  \end{equation*}    
Then we have
 \begin{equation*}
      \partial \f r(\xi) = \partial \f s(\f t(\xi)) \partial \f t(\xi),
 \end{equation*}
and consequently, for every function $v:S \rightarrow \R,$
 \begin{equation*}
 v\circ \f r:\Omega \rightarrow \R \quad \text{and} \quad v\circ \f s:\Omega_s\rightarrow \R,
 \end{equation*}
we have
   \begin{equation}
   \label{eq:reparam}
  \nabla_g v(\f r(\xi))\; \partial \f r(\xi) = \nabla (v \circ \f r)(\xi) \text{ and  } 
   \nabla_g v(\f s(\zeta))\; \partial \f s(\zeta) = \nabla (v \circ \f s )(\zeta) .
   \end{equation}
Using chain rule on both sides of the former equation of \eqref{eq:reparam}, then we get
\begin{equation*}
       \nabla_g v(\f s(\zeta))\; \partial \f s(\f t(\xi)) \partial\f t (\xi) = 
       \partial (v \circ \f s(\f t(\xi))) \partial \f t(\xi) \Rightarrow  \nabla_g v(\f s(\zeta))\; \partial \f s(\f t(\xi)) =\partial (v \circ \f s(\f t(\xi))) ,
\end{equation*}
which gives the latter equation in \eqref{eq:reparam} since $\partial \f t(\xi)$ is non-degenerate.
Using the same process but consider $\f t^{-1}:\Omega_s \rightarrow \Omega$, we can show the reverse implication.
Thus, we have shown that any two arbitrary parameterizations $\f r$ and $\f s$ lead to the same gradient values at same positions.
\end{proof}

\newpage

\begin{thebibliography}{10}

\bibitem{AinsworthOden2000}
{\sc M.~Ainsworth and J.~T. Oden}, {\em A posteriori error estimation in finite
  element analysis}, Pure and Applied Mathematics (New York),
  Wiley-Interscience [John Wiley \& Sons], New York, 2000.

\bibitem{Aubin1982}
{\sc T.~Aubin}, {\em Best constants in the {S}obolev imbedding theorem: the
  {Y}amabe problem}, in Seminar on {D}ifferential {G}eometry, vol.~102 of Ann.
  of Math. Stud., Princeton Univ. Press, Princeton, N.J., 1982, pp.~173--184.

\bibitem{Babuska2001}
{\sc I.~Babu{\v{s}}ka and T.~Strouboulis}, {\em The finite element method and
  its reliability}, Numerical Mathematics and Scientific Computation, The
  Clarendon Press, Oxford University Press, New York, 2001.

\bibitem{BankXu2003}
{\sc R.~E. Bank and J.~Xu}, {\em Asymptotically exact a posteriori error
  estimators. {I}. {G}rids with superconvergence}, SIAM J. Numer. Anal., 41
  (2003), pp.~2294--2312 (electronic).

\bibitem{CamachoDemlow2015}
{\sc F.~Camacho and A.~Demlow}, {\em {$L_2$} and pointwise a posteriori error
  estimates for {FEM} for elliptic {PDE}s on surfaces}, IMA J. Numer. Anal., 35
  (2015), pp.~1199--1227.

\bibitem{ChenGuoZhangZou2017}
{\sc H.~Chen, H.~Guo, Z.~Zhang, and Q.~Zou}, {\em A {$C^0$} linear finite
  element method for two fourth-order eigenvalue problems}, IMA J. Numer.
  Anal., 37 (2017), pp.~2120--2138.

\bibitem{Chen2008}
{\sc L.~Chen}, {\em Short implementation of bisection in {MATLAB}}, in Recent
  advances in computational sciences, World Sci. Publ., Hackensack, NJ, 2008,
  pp.~318--332.

\bibitem{chen2009ifem}
\leavevmode\vrule height 2pt depth -1.6pt width 23pt, {\em {$i$FEM}: an
  innovative finite element methods package in {MATLAB}}, University of
  California at Irvine,  (2009).

\bibitem{Chernyshenk02015}
{\sc A.~Y. Chernyshenko and M.~A. Olshanskii}, {\em An adaptive octree finite
  element method for {PDE}s posed on surfaces}, Comput. Methods Appl. Mech.
  Engrg., 291 (2015), pp.~146--172.

\bibitem{DednerMadhavanStinner2016}
{\sc A.~Dedner and P.~Madhavan}, {\em Adaptive discontinuous {G}alerkin methods
  on surfaces}, Numer. Math., 132 (2016), pp.~369--398.

\bibitem{DednerMadhavanStinner2013}
{\sc A.~Dedner, P.~Madhavan, and B.~Stinner}, {\em Analysis of the
  discontinuous {G}alerkin method for elliptic problems on surfaces}, IMA J.
  Numer. Anal., 33 (2013), pp.~952--973.

\bibitem{Demlow2009}
{\sc A.~Demlow}, {\em Higher-order finite element methods and pointwise error
  estimates for elliptic problems on surfaces}, SIAM J. Numer. Anal., 47
  (2009), pp.~805--827.

\bibitem{DemlowDziuk2007}
{\sc A.~Demlow and G.~Dziuk}, {\em An adaptive finite element method for the
  {L}aplace-{B}eltrami operator on implicitly defined surfaces}, SIAM J. Numer.
  Anal., 45 (2007), pp.~421--442 (electronic).

\bibitem{Demlow2012}
{\sc A.~Demlow and M.~A. Olshanskii}, {\em An adaptive surface finite element
  method based on volume meshes}, SIAM J. Numer. Anal., 50 (2012),
  pp.~1624--1647.

\bibitem{doCarmo1992}
{\sc M.~P. do~Carmo}, {\em Riemannian geometry}, Mathematics: Theory \&
  Applications, Birkh\"auser Boston, Inc., Boston, MA, 1992.
\newblock Translated from the second Portuguese edition by Francis Flaherty.

\bibitem{DonJueSchTak17}
{\sc G.~Dong, B.~J{\"u}ttler, O.~Scherzer, and T.~Takacs}, {\em Convergence of
  tikhonov regularization for solving ill--posed operator equations with
  solutions defined on surfaces}, Inverse Probl. Imaging, 11 (2017), pp.~221 --
  246.

\bibitem{Dorfler1996}
{\sc W.~D{\"o}rfler}, {\em A convergent adaptive algorithm for {P}oisson's
  equation}, SIAM J. Numer. Anal., 33 (1996), pp.~1106--1124.

\bibitem{DuJu2005}
{\sc Q.~Du and L.~Ju}, {\em Finite volume methods on spheres and spherical
  centroidal {V}oronoi meshes}, SIAM J. Numer. Anal., 43 (2005), pp.~1673--1692
  (electronic).

\bibitem{Dziuk1988}
{\sc G.~Dziuk}, {\em Finite elements for the {B}eltrami operator on arbitrary
  surfaces}, in Partial differential equations and calculus of variations,
  vol.~1357 of Lecture Notes in Math., Springer, Berlin, 1988, pp.~142--155.

\bibitem{DziukElliott2013}
{\sc G.~Dziuk and C.~M. Elliott}, {\em Finite element methods for surface
  {PDE}s}, Acta Numer., 22 (2013), pp.~289--396.

\bibitem{GrandeReusken2016}
{\sc J.~Grande and A.~Reusken}, {\em A higher order finite element method for
  partial differential equations on surfaces}, SIAM J. Numer. Anal., 54 (2016),
  pp.~388--414.

\bibitem{GuoZhang2015}
{\sc H.~Guo and Z.~Zhang}, {\em Gradient recovery for the {C}rouzeix-{R}aviart
  element}, J. Sci. Comput., 64 (2015), pp.~456--476.

\bibitem{GuoZhangZhao2016}
{\sc H.~Guo, Z.~Zhang, and R.~Zhao}, {\em Hessian recovery for finite element
  methods}, Math. Comp., 86 (2017), pp.~1671--1692.

\bibitem{GuoZhangZhao2017}
\leavevmode\vrule height 2pt depth -1.6pt width 23pt, {\em Superconvergent
  two-grid methods for elliptic eigenvalue problems}, J. Sci. Comput., 70
  (2017), pp.~125--148.

\bibitem{GuoZhangZou2017}
{\sc H.~Guo, Z.~Zhang, and Q~Zou}, {\em A {$C^0$} {L}inear {F}inite {E}lement
  {M}ethod for {B}iharmonic {P}roblems}, J. Sci. Comput., 74 (2018),
  pp.~1397--1422.

\bibitem{GZZ2018b}
{\sc H.~Guo, Z.~Zhang, and Q.~Zou}, {\em A {$C^0$} linear finite element method
  for sixth order elliptic equations}, arXiv:1804.03793 [math.NA], 2018.

\bibitem{Hebey1999}
{\sc E.~Hebey}, {\em Nonlinear analysis on manifolds: {S}obolev spaces and
  inequalities}, vol.~5 of Courant Lecture Notes in Mathematics, New York
  University, Courant Institute of Mathematical Sciences, New York; American
  Mathematical Society, Providence, RI, 1999.

\bibitem{Lakhany2000}
{\sc A.~M. Lakhany, I.~Marek, and J.~R. Whiteman}, {\em Superconvergence
  results on mildly structured triangulations}, Comput. Methods Appl. Mech.
  Engrg., 189 (2000), pp.~1--75.

\bibitem{Lee2013}
{\sc J.~M. Lee}, {\em Riemannian manifolds}, vol.~176 of Graduate Texts in
  Mathematics, Springer-Verlag, New York, 1997.
\newblock An introduction to curvature.

\bibitem{NagaZhang2004}
{\sc A.~Naga and Z.~Zhang}, {\em A posteriori error estimates based on the
  polynomial preserving recovery}, SIAM J. Numer. Anal., 42 (2004),
  pp.~1780--1800 (electronic).

\bibitem{NagaZhang2012}
\leavevmode\vrule height 2pt depth -1.6pt width 23pt, {\em Function value
  recovery and its application in eigenvalue problems}, SIAM J. Numer. Anal.,
  50 (2012), pp.~272--286.

\bibitem{NagaZhangZhou2006}
{\sc A.~Naga, Z.~Zhang, and A.~Zhou}, {\em Enhancing eigenvalue approximation
  by gradient recovery}, SIAM J. Sci. Comput., 28 (2006), pp.~1289--1300.

\bibitem{OlshankiiReusken2010}
{\sc M.~A. Olshanskii and A.~Reusken}, {\em A finite element method for surface
  {PDE}s: matrix properties}, Numer. Math., 114 (2010), pp.~491--520.

\bibitem{OlshanskiiReuskenGrande2009}
{\sc M.~A. Olshanskii, A.~Reusken, and J.~Grande}, {\em A finite element method
  for elliptic equations on surfaces}, SIAM J. Numer. Anal., 47 (2009),
  pp.~3339--3358.

\bibitem{OlshanskiiSafin2016}
{\sc M.~A. Olshanskii and D.~Safin}, {\em A narrow-band unfitted finite element
  method for elliptic {PDE}s posed on surfaces}, Math. Comp., 85 (2016),
  pp.~1549--1570.

\bibitem{cgal}
{\sc L.~Rineau and M.~Yvinec}, {\em {3D} surface mesh generation}, in {CGAL}
  User and Reference Manual, {CGAL Editorial Board}, {4.9}~ed., 2016.

\bibitem{WeiChenHuang2010}
{\sc H.~Wei, L.~Chen, and Y.~Huang}, {\em Superconvergence and gradient
  recovery of linear finite elements for the {L}aplace-{B}eltrami operator on
  general surfaces}, SIAM J. Numer. Anal., 48 (2010), pp.~1920--1943.

\bibitem{XuZhang2004}
{\sc J.~Xu and Z.~Zhang}, {\em Analysis of recovery type a posteriori error
  estimators for mildly structured grids}, Math. Comp., 73 (2004),
  pp.~1139--1152 (electronic).

\bibitem{XGZ2019}
{\sc M.~Xu, H.~Guo, and Q.~Zou}, {\em Hessian recovery based finite element
  methods for the cahn-hilliard equation}, J. Comput. Phys.,  (2019).

\bibitem{ZhangNaga2005}
{\sc Z.~Zhang and A.~Naga}, {\em A new finite element gradient recovery method:
  superconvergence property}, SIAM J. Sci. Comput., 26 (2005), pp.~1192--1213
  (electronic).

\bibitem{ZZ1992}
{\sc O.~C. Zienkiewicz and J.~Z. Zhu}, {\em The superconvergent patch recovery
  and a posteriori error estimates. {I}. {T}he recovery technique}, Internat.
  J. Numer. Methods Engrg., 33 (1992), pp.~1331--1364.

\bibitem{ZZ1992b}
\leavevmode\vrule height 2pt depth -1.6pt width 23pt, {\em The superconvergent
  patch recovery and a posteriori error estimates. {II}. {E}rror estimates and
  adaptivity}, Internat. J. Numer. Methods Engrg., 33 (1992), pp.~1365--1382.

\end{thebibliography}
%

\end{document}